\documentclass[peprint,3p,times]{elsarticle}
\usepackage[utf8]{inputenc}
\usepackage[normalem]{ulem}
\usepackage[size=tiny]{todonotes}
\usepackage{tikz}
\usepackage{colortbl}
\usepackage{soul}
\usetikzlibrary{positioning}

\newcommand{\stkout}[1]{\ifmmode\text{\sout{\ensuremath{#1}}}\else\sout{#1}\fi}

\def\bff{{\bm f}}
\def\hbff{\widehat{\bm f}}
\def\tbff{\widetilde{\bff}}

\def\bphi{\bm \phi}
\def\bPsi{\bm \Psi}
\def\bpsi{\bm \psi}
\def\beta{\bm \eta}

\def\bE{{\bm E}}
\def\tbE{\widetilde{\bE}}
\def\hbphi{\widehat{\bm \phi}}
\def\tbphi{\widetilde{\bphi}}

\def\bg{\bm g}

\def\bu{\bm u}
\def\bp{\bm p}
\def\bv{\bm v}

\def\wbG{\widehat{{\bm G}}}
\def\bN{{\bm N}}
\def\N{{\mathbb{N}}}

\def\wbA{\widehat{\bm{A}}}

\def\bDv{{\bm D}_v}
\def\bDcx{{\bm D}_{c,x}}
\def\bDcv{{\bm D}_{c,v}}

\def\bIx{{\bm I_x}}
\def\bIv{{\bm I_v}}
\def\bF{{\bm F}}
\def\bU{{\bm U}}
\def\bSigma{{\bm \Sigma}}
\def\bV{{\bm V}}
\def\tdiag{\text{diag}}
\def\tos{\text{os}}
\def\real{\mathbb{R}}

\def\bmG{\widehat{\boldsymbol{\mathcal{G}}}}

\def\vT{v_T}
\def\mD{\mathcal{D}}
\def\Tinterval{T_{\textup{interval}}}
\def\epf{\epsilon_{f,t_f}}
\def\bptrain{\bp_{\text{train}}}
\def\bptest{\bp_{\text{test}}}

\def\mF{\mathcal{F}}
\def\mS{\mathcal{S}}
\def\mG{\mathcal{G}}
\def\wbg{\widehat{\bm{g}}}
\def\mT{\mathcal{T}}
\def\bfk{{\bm k}}

\newcommand{\tento}[1]{10^{#1}}

\newcommand{\arbfunc}{\Theta}
\newcommand{\dash}{--}
\newcommand{\indicator}{\Xi}

\usepackage[english]{babel}

\usepackage{amsmath,bm,mathtools}
\usepackage{amssymb}
\usepackage{graphicx}
\usepackage[colorlinks=true, allcolors=blue]{hyperref}
\usepackage{subcaption}
\usepackage{comment}

\newtheorem{remark}{Remark}

\begin{document}

\begin{frontmatter}

%\title{Accelerating Kinetic Simulations of Electrostatic Plasmas with Reduced-Order Modeling\tnoteref{t1}}
\title{
Local Reduced-Order Modeling for Electrostatic Plasmas by Physics-Informed Solution Manifold Decomposition\tnoteref{t1}}
\tnotetext[t1]{Lawrence Livermore National Laboratory is operated by Lawrence Livermore National Security, LLC, for the U.S. Department of Energy, National Nuclear Security Administration under Contract DE-AC52-07NA27344}

\author[1]{Ping-Hsuan Tsai \corref{cor1}}
\ead{pinghsuan@vt.edu}

\author[2]{Seung Whan Chung}
\ead{chung28@llnl.gov}

\author[2]{Debojyoti Ghosh}
\ead{ghosh5@llnl.gov}

\author[2]{John Loffeld}
\ead{loffeld1@llnl.gov}

\author[2]{Youngsoo Choi}
\ead{choi15@llnl.gov}

\author[2]{Jonathan L. Belof}
\ead{belof1@llnl.gov}

\cortext[cor1]{Corresponding author}

\affiliation[1]{organization={Department of Mathematics, Virginia Tech},%Department and Organization
            city={Blacksburg},
            postcode={24061}, 
            state={VA},
            country={United States}}
            
\affiliation[2]{organization={ Lawrence Livermore National Laboratory},
city={Livermore}, postcode={94550}, state={CA}, country={United States}}

\begin{abstract}
  Despite advancements in high-performance computing and modern numerical
    algorithms, computational cost remains prohibitive for multi-query kinetic plasma
    simulations. In this work, we develop data-driven reduced-order models (ROMs)
    for collisionless electrostatic plasma dynamics, based on the kinetic
    Vlasov-Poisson equation. Our ROM approach projects the equation onto a
    linear subspace defined by the proper orthogonal decomposition (POD)
    modes. We introduce an efficient tensorial method to update the nonlinear
    term using a precomputed third-order tensor. We capture multiscale behavior
    with a minimal number of POD modes by decomposing the solution manifold into multiple
    time windows and creating temporally local
    ROMs. We consider two strategies for decomposition: one based on the physical time and the other based on the electric field energy.
    Applied to the 1D1V Vlasov\dash Poisson simulations, that is, prescribed E-field, Landau damping, and two-stream instability, 
    we demonstrate that our ROMs accurately capture the total energy of the system both for parametric and time extrapolation cases.
    The temporally local ROMs are more efficient and accurate than the single ROM. 
    In addition, in the two-stream instability case, we show that the energy-windowing reduced-order model (EW-ROM) is more efficient and accurate than the time-windowing reduced-order model (TW-ROM). With the tensorial approach, EW-ROM solves the equation approximately $90$ times faster than Eulerian simulations while maintaining a maximum relative error of $7.5\%$ for the training data and $11\%$ for the testing data.
\end{abstract}


\end{frontmatter}

\section{Introduction}

Collisionless electrostatic plasma kinetics is governed by the Vlasov\dash Poisson equation, which describes the evolution of charged particle distribution under the self-consistent electrostatic field.
Due to high dimensionality, scale disparities, and nonlinearity, solving these equations is computationally challenging. Lagrangian particle-in-cell (PIC) methods~\cite{birdsall2018plasma} are often used, where particles are advanced along the characteristic curves of the Vlasov equation~\cite{verboncoeur2005particle}. 
Alternatively, grid-based Eulerian methods \cite{hypar,dorf2012progress,lesser2022loki,vogman2014dory} allow the use of advanced high-order numerical algorithms for partial differential equations (PDEs)~\cite{hittinger2013block,sonnendrucker1999semi,einkemmer2016high}.
In this work, we consider the 1D1V Vlasov\dash Poisson equation in the parametric setting. In particular, we study the effects of initial conditions  
characterized by the perturbation amplitude and the thermal velocity. 
While high-performance computing capabilities and scalable algorithms enable high-fidelity kinetic simulations, parametric studies remain intractable due to requiring many forward simulations. 

Reduced order modeling can be a promising alternative for such multi-query applications,
by providing a surrogate model that solves the equation with moderate accuracy but considerable speed-up compared to the corresponding full-order model (FOM).
Reduced-order models (ROMs) can be divided into two categories, namely, intrusive and nonintrusive ROMs. 
Intrusive ROMs are constructed using knowledge of underlying governing equations, numerical schemes, and solution data, whereas nonintrusive ROMs are constructed using solution data only \cite{rozza2022advanced,fries2022lasdi,he2023glasdi,bonneville2024gplasdi,tran2024weak,bonneville2024comprehensive,kim2024gappy}. 
We consider projection-based reduced-order models (ROMs), which is one of the intrusive approaches. In particular, the reduced basis vectors are extracted by performing proper orthogonal decomposition (POD) on the snapshot data of the FOM simulations, and the governing equations are projected onto the reduced space. These approaches take advantage of both the \textit{known governing equations} and the \textit{solution data} generated from the corresponding FOM simulations to form linear subspace reduced order models (LS-ROM).
In spite of the success of the classical LS-ROMs in many applications, e.g., %have been successfully applied to many physical simulations, e.g., 
fluid dynamics
\cite{choi2019space,carlberg2018conservative,kim2022fast,lauzon2024s,diaz2024fast,kaneko2020towards,tsai2022parametric},
nonlinear diffusion \cite{hoang2021domain}, Boltzmann transport
\cite{choi2021space}, design optimization
\cite{choi2020gradient,mcbane2021component,mcbane2022stress,choi2019accelerating},
these approaches are limited to the assumption that the intrinsic solution space
falls into a subspace with a small dimension, that is, the solution space with a
Kolmogorov $n$-width decays fast. This assumption is violated in
advection-dominated and advection problems
\cite{tsai2023accelerating,tsai2024time,kim2022fast,zanardi2024scalable,cheung2023local}, which prevents reduced-order modeling from being practical.
Recently, there have been many attempts to develop efficient ROMs for
advection-dominated problems, and these can be divided into two categories.
The first category replaces the linear subspace solution representation with a
nonlinear manifold \cite{kim2022fast,RAZAVI2025113576}. 
The second category enhances the solution representability of the linear subspace by introducing special treatments and adaptive schemes.
We note that although the nonlinear manifold approach is powerful for
advection-dominated problems, the training of the neural network could be intricate  
compared to linear subspace methods, especially in large-scale problems.

In this paper, we consider the idea of building small and accurate reduced-order models by decomposing the solution manifold into submanifolds, which belongs to the second category.
These reduced-order models are local in the sense that each of them is valid only over a certain subset of the parameter-time domain. The appropriate local reduced-order model is chosen based on representative quantities of the current state of the system. The concept of a local reduced-order model is introduced in \cite{washabaugh2012nonlinear,amsallem2012nonlinear}, where unsupervised clustering is used for the solution manifold decomposition. In \cite{copeland2022reduced}, the authors used physical time as a solution manifold indicator, and the resulting time-windowing ROM was able to achieve good speed-up with accurate approximated solutions in various benchmark problems, including Sedov blast, Gresho vortices, Taylor\dash Green vortices, and triple-point problems. In \cite{cheung2023local}, the authors considered the penetration speed of a fluid interface as the physical indicator for solution decomposition, and the resulting local ROMs yield better approximations than the time-windowing ROMs in the Rayleigh\dash Taylor instability.
In \cite{cheung2024data}, the authors considered reduced-order models with dynamic mode decomposition (DMD) and solution manifold decomposition and demonstrated the time-windowing DMD is able to reproduce the 
complex pore collapse processes and predict to a certain extent. 
In \cite{shimizu2021windowed,parish2021windowed}, the authors consider the time-windowing approach with space-time least-squares Petrov\dash Galerkin method for model reduction of dynamical systems.
In alignment with the concept of decomposing the solution manifold into submanifolds, in \cite{hesthaven2023adaptive}, the authors developed projection-based reduced-order models for the parametric 1D1V Vlasov--Poisson equations. However, we note that these ROMs 
were developed for the semi-discrete Hamiltonian system resulting from a geometric particle-in-cell approximation. In addition, the ROMs are based on the dynamical low-rank approximation where the reduced space evolves over time. This approach differs fundamentally from the one we are proposing.

Our contribution in this paper is to propose an efficient model order reduction scheme for accelerating the kinetic simulations of electrostatic plasmas with varying perturbation amplitude and thermal velocity in the initial conditions.
Similarly to \cite{cheung2023local}, the idea is to construct temporally local projection-based ROMs that are small but accurate within a short period in advection-dominated problems to achieve a good speed-up and solution accuracy.
Despite the temporally local ROMs being small, the primary computational expense arises from the nonlinear hyperbolic term evaluation in the high-dimensional space. To improve efficiency, we introduce a tensorial approach for updating this term using a precomputed third-order tensor. 
In this work, we consider the physical time and the electric field energy as the indicators for the solution manifold decomposition and compare the performance in speed-up and solution accuracy.
Our numerical experiments demonstrate the efficacy of the temporally local ROMs compared to the single linear subspace ROM. In addition, the electric field energy is a good alternative indicator to resolve the deficiency of degenerating solution accuracy with respect to perturbation amplitude and thermal velocity in the time-windowing approach.

The rest of the paper is organized as follows. In Section~\ref{sec:FOM}, we introduce the governing equations and the numerical discretization that will be used as the full-order model. Next, a projection-based ROM and the tensorial approach, which is used to efficiently update the nonlinear term, are described in Section~\ref{sec:ROMs}. In Section~\ref{sec:temporal_ROM}, time-windowing and energy-windowing reduced-order models are presented. Numerical results are presented in Section~\ref{sec:numerical}. Finally, conclusions are given in Section~\ref{sec:conclusions}.

\section{Full Order Model (FOM)}
\label{sec:FOM}
We consider the 1D1V Vlasov\dash Poisson equation,
\begin{align} 
    & \partial_t f + v \partial_x f + E\partial_v f= 0 \label{eq:1D--1V-vlasov} \\
    & \partial^2_x \phi = \displaystyle\int f~dv 
    - 
    \frac{1} {\int dx}
    \int \int f~dv~dx
    %\quad \text{and}  \quad
    , \label{eq:1D--1V-poisson}
\end{align}
where $f(x,v,t)$ is the probability distribution function of electrons, $E(x, t)$ is the electric field, and $\phi(x, t)$ is the self-consistent electrostatic potential with $E = -\partial_x \phi$. $x$ and $v$ are the spatial coordinate and the velocity, respectively. 
The second term of the right-hand side in (\ref{eq:1D--1V-poisson}) corresponds to ion charge density, which is assumed to be uniform and constant in space and time due to its heavy mass compared to electron.
The boundary conditions are periodic in the $x$-direction and homogeneous Dirichlet in the $v$-direction. In order to close the system, an initial condition needs to be imposed. In this work, the initial condition is considered in the parametric setting by regarding the perturbation amplitude $\alpha$ and thermal velocity $v_T$ as problem parameters, that is, $f(x,v,0;\alpha,v_T)$ in which $(\alpha, v_T)$ lies in a parametric domain $\mathcal{D}\subset \mathbb{R}^2$.

The governing equation (\ref{eq:1D--1V-vlasov}) can be rewritten as a hyperbolic PDE
\begin{align}
    \frac{\partial {f}}{\partial t} 
    +
    \frac{\partial}{\partial x}  { F}\left({ f}\right)
    +
    \frac{\partial}{\partial v} { G}\left({ f}\right)
    = 0,
    \ \ \ \ 
    { F}\left({ f}\right) \equiv v { f},\ \ \ \ 
    { G}\left({ f}\right) \equiv E { f} \label{eq:hyp_pde}.
\end{align}
The spatial derivatives in (\ref{eq:hyp_pde}) are discretized using the conservative finite-difference formulation~\cite{jiangshu1996}, which is summarized below for an arbitrary function $\arbfunc\left(f\right)$ in 1D.
The $x$-derivative $\partial_x \arbfunc$ at grid point $x_i$ is written as
\begin{align}
    \left.\frac{\partial \arbfunc}{\partial x}\right|_i
    =
    \frac{1}{\Delta x}
    \left(\hat{\arbfunc}_{i+\frac{1}{2}} - \hat{\arbfunc}_{i-\frac{1}{2}}\right)
    +
    \mathcal{O}\left(\Delta x^p\right), \label{eq:cfd}
\end{align}
where the locations $i\pm1/2$ denote the cell interfaces of grid cell $i$ 
, $p$ is the order the spatial discretization, $\hat{\arbfunc}$ is the numerical approximation to the primitive of $\arbfunc$, and $\Delta x$ is the grid spacing. The upwind flux at the interface, $\hat{\arbfunc}_{i\pm\frac{1}{2}}$, is computed as follows:
\begin{align}
    \hat{\arbfunc}_{i+\frac{1}{2}} = 
    \begin{cases}
        \hat{\arbfunc}^L_{i+\frac{1}{2}} & \arbfunc'_{j} > 0, \arbfunc'_{j+1} > 0 \\
        \hat{\arbfunc}^R_{i+\frac{1}{2}} & \arbfunc'_{j} < 0, \arbfunc'_{j+1} < 0 \\
        \frac{1}{2}\left[\hat{\arbfunc}^L_{i+\frac{1}{2}} + \hat{\arbfunc}^R_{i+\frac{1}{2}} - \left(\max\limits_{j,j+1}{\arbfunc'}\right) \left(\hat{f}^R_{i+\frac{1}{2}} - \hat{f}^L_{i+\frac{1}{2}}\right)\right] & {\rm otherwise},
    \end{cases}
    \label{eq:upwindflux}
\end{align}
where $\arbfunc'$ denotes the derivative of $\arbfunc\left(f\right)$ with respect to $f$. The scalar terms $\arbfunc^{L,R}_{i+\frac{1}{2}}$ and $\hat{f}^{L,R}_{i+\frac{1}{2}}$ represent the left- and right-biased approximations to $\arbfunc\left(f\right)$ and $f$ at $i+\frac{1}{2}$, respectively. They are computed with the 5th-order weighted essentially nonoscillatory (WENO5) scheme ~\cite{jiangshu1996} and its linear counterpart (UPWIND5). The algorithm summarized above is applied to $\partial_x {F}\left({f}\right)$ and $\partial_v { F}\left({f}\right)$ in~(\ref{eq:hyp_pde}). The resulting spatially discretized ordinary differential equation (ODE) is 
\begin{align}
    \begin{cases}
    \displaystyle\frac{d{\bff}}{dt} 
    +
    {\bN}({\bff}) = {\bf 0} \\ 
    \bphi = \mathcal{F}^{-1} \left(-\bfk^{-2} \mathcal{F}\left({\bg}(\bff)\right)\right),
    \end{cases}
    \label{eq:ode}
\end{align}
where ${\bN}(\bff)$ is a vector represents the hyperbolic term approximated by the WENO5/UPWIND5 scheme, and $\bg(\bff)$ is the discretization of the right-hand side in (\ref{eq:1D--1V-poisson}).
The solution vector ${\bff \in \real^{N_f=n_x\times n_v}}$ comprises the distribution function $f\left(x_i, v_j, t\right)$ at each grid point $x_i, v_j$, where $ 0 \le i \le n_x, 0 \le j \le n_v$ and the solution vector $\bphi \in \real^{N_\phi = n_x}$ comprises the potential function $\phi(x_i,t)$ at each grid point $x_i$, where $ 0 \le i \le n_x$.
The ODEs are solved using the classical 4th-order, 4-stage explicit Runge\dash Kutta method. In this paper, we only consider simulations that are periodic in space; thus,~(\ref{eq:1D--1V-poisson}) is solved exactly using the Fast Fourier Transform (FFT), and $\mathcal{F}$ and $\mathcal{F}^{-1}$ are the Fourier transform and inverse Fourier transform. The FOM is implemented in {\tt HyPar}~\cite{hypar}, an MPI-based C/C++ code for hyperbolic and parabolic PDEs, and the {\tt FFTW} library~\cite{FFTW05} is used for the FFT-based Poisson solve.

\section{Reduced Order Models (ROMs)}
\label{sec:ROMs}

In this section, we present the details of our projection-based reduced order
model for the semi-discrete 1D1V Vlasov\dash Poisson equation
(\ref{eq:ode}).
The reduced order model is constructed in the offline phase and deployed in the
online phase. In what follows, we first discuss the construction of the reduced
basis functions and then move on to discuss all of the essential ingredients of
the reduced order model.

\subsection{Solution bases construction} \label{sec:solution_basis_construction}

In the reduced order model, we restrict the solution spaces of the probability distribution, potential, and electric fields to subspaces spanned by the reduced basis functions. %for the probability distribution function $f$. 
That is, the subspace for the probability distribution is defined as
\begin{equation}
   \mS_f \equiv \text{Span}\{\bpsi^i_f\}^{n_f}_{i=1} \subseteq \real^{N_f} ,
\end{equation}
with dim$(\mathcal{S}_f) = n_f \ll N_f$. The  
subspaces for the potential and electric fields are defined as 
\begin{equation}
    \mathcal{S}_\phi \equiv \text{Span}\{\bpsi^i_{\phi}\}^{n_\phi}_{i=1} \subset \mathbb{R}^{N_\phi}, \quad
    \mathcal{S}_E \equiv \text{Span}\{\bpsi^i_{E}\}^{n_\phi}_{i=1} \subset \mathbb{R}^{N_\phi},
\end{equation}
with $\text{dim}(\mathcal{S}_{\phi}) = \text{dim}(\mathcal{S}_E) = n_\phi \ll N_\phi$. 
Using these subspaces, the discrete fields $\bff$, $\bphi$, and $\bE$ are approximated in trial subspaces:
\begin{equation}
    \bff \approx \tbff = \bff_{\tos} + \bPsi_f \hbff, %\in \bff_\tos + \mathcal{S}_f, 
    \quad \bphi \approx \tbphi = \bphi_\tos + \bPsi_\phi \hbphi, %\in \bphi_\tos + \mathcal{S}_\phi,
    \quad \bE \approx \tbE = \bE_\tos + \bPsi_E \hbphi, %\in \bphi_\tos + \mathcal{S}_\phi, 
    \label{eq:reduced_app}
\end{equation}
where $\bff_\tos \in \mathbb{R}^{N_f}$, $\bphi_\tos \in \mathbb{R}^{N_\phi}$, and $\bE_\tos \in \mathbb{R}^{N_\phi}$ denote the prescribed offset vector for the distribution, potential, and electric fields, respectively. The orthonormal basis matrices $\bPsi_f \in \mathbb{R}^{N_f \times n_f}$, $\bPsi_\phi \in \mathbb{R}^{N_\phi \times n_\phi}$, and $\bPsi_E \in \mathbb{R}^{N_\phi \times n_\phi}$ %for the distribution, potential, and electric fields are defined as 
are defined as
\begin{equation}
    \bPsi_f \equiv \begin{bmatrix}
    \bpsi^1_f & \bpsi^2_f & \cdots & \bpsi^{n_f}_{f}
    \end{bmatrix}, 
    \quad \bPsi_\phi \equiv 
    \begin{bmatrix}
    \bpsi^1_\phi & \bpsi^2_\phi & \cdots & \bpsi^{n_\phi}_{\phi}\end{bmatrix},
    \quad \bPsi_E \equiv 
    \begin{bmatrix}
    \bpsi^1_E & \bpsi^2_E & \cdots & \bpsi^{n_\phi}_{E}\end{bmatrix},
\end{equation}
and $\hbff:[0,t_f] \rightarrow \mathbb{R}^{n_f}$ and $\hbphi:[0,t_f] \rightarrow \mathbb{R}^{n_\phi}$ denote the time-dependent generalized coordinates of the distribution and the potential, respectively.

We consider proper orthogonal decomposition (POD) \cite{berkooz1993proper} to
construct the reduced basis matrix $\bPsi_f$ for the distribution from a
truncated singular
value decomposition (SVD) approximation to the FOM distribution snapshot matrix.
We collect the solution data by running FOM simulations on a set of problem
parameters, namely $\{\bp_k = (\alpha_k, v_{T,k})\}^{n_p}_{k=1}$. In this paper,
we consider the same final time $t_f$ and the number of time steps $N_t$ for all
the problem parameters in the FOM simulation. A solution snapshot matrix is
formed by assembling all the FOM solution data, that is\footnote{We note that
the boldface $\bF$ denotes snapshot matrix and is different from $F$ in
(\ref{eq:hyp_pde}), which is denoted as the flux.},
\begin{equation}
    \bF \equiv \begin{bmatrix}
        \bff_1(\bp_1)& \cdots \bff_{N_t}(\bp_{1}) & \cdots & \bff_{1}(\bp_{n_p}) & \cdots & \bff_{N_t}(\bp_{n_p})
    \end{bmatrix} \in \mathbb{R}^{N_f\times K}, \label{eq:full_snapshot_matrix}
\end{equation}
where $\bff_n(\bp_{k})$ is the distribution state at the $n$th time step with problem parameter $\bp_k$ for $n\in \N(N_t)$ computed from the FOM simulation and $K=\sum^{n_p}_{k=1}N_t$.  Then, POD computes its thin SVD,
\begin{equation*}
    \bF = \bU \bSigma \bV^T,
\end{equation*}
where $\bU \in \real^{N_f \times K}$ and $\bV = \real^{K \times K}$ are orthogonal matrices, and $\bSigma \in \real^{K \times K}$ is the diagonal singular value matrix. The basis matrix $\bPsi_f = \begin{bmatrix}
    \bu_1 ~\ldots~ \bu_{n_f}
\end{bmatrix}$ is set to be the leading $n_f$ columns of $\bU$. The basis size, $n_f$, is determined by the energy criteria, that is, we find the minimum $n_f \in \mathbb{N}(K)$ such that the following condition is satisfied:
\begin{equation}
    \delta_\sigma \ge 1 -\frac{\sum^{n_f}_{i=1}\sigma_i}{\sum^{K}_{i=1}\sigma_i}, \label{eq:energy_criteria}
\end{equation}
where $\sigma_i$ is the $i$-th largest singular value in the singular matrix, $\bSigma$, and $\delta_\sigma \in [0,1]$ is the energy missing ratio. %denotes a threshold. 
The POD procedure seeks the optimal $n_f$-dimensional subspace by minimizing
$\|\bF - \bPsi_f \bPsi^T_f \bF\|^2_F$ in the Frobenius norm over all $\bPsi_f
\in \mathbb{R}^{N_f \times K}$ with orthonormal columns. One could also seek the optimal subspace in terms of the $L^2$ and $H^1_0$ norms \cite{kaneko2020towards,tsai2022parametric,fick2018stabilized}.

We note that there are two approaches to construct the potential reduced basis matrix $\bPsi_\phi$.
In the first approach, we collect the solution data on a set of problem
parameters $\{\bp_k = (\alpha_k, v_{T,k})\}^{n_p}_{k=1}$ and form a solution snapshot matrix, 
\begin{equation}
    \bF_\phi \equiv \begin{bmatrix}
        \bphi_1(\bp_1) & \cdots & \bphi_{N_t}(\bp_{1}) & \cdots & \bphi_{1}(\bp_{n_p}) & \cdots & \bphi_{N_t}(\bp_{n_p})
    \end{bmatrix} \in \mathbb{R}^{N_\phi\times K},
\end{equation}
where $\bphi_n(\bp_{k})$ is the potential state at the $n$th time step with problem parameter $\bp_k$ for $n\in \N(N_t)$ computed from the FOM simulation and $K=\sum^{n_p}_{k=1}N_t$.  Then, POD computes the thin SVD of $\bF_\phi$, and the basis matrix $\bPsi_\phi$ is set to be the leading $n_\phi$ columns of the left singular matrix. The basis size, $n_\phi$, is determined similarly as in the distribution using (\ref{eq:energy_criteria}).

In the second approach, the {prescribed offset vector $\bphi_\tos$} and the reduced basis matrix $\bPsi_\phi$ are constructed by solving Poisson problems using FFT with the prescribed offset vector $\bff_\tos$ and the distribution reduced basis functions $\bpsi^i_f$, that is, 
\begin{equation}
    \bphi_\tos = -\mF^{-1}(\bfk^{-2} \mF(\bg(\bff_\tos))), \quad \bpsi^i_{\phi} = -\mF^{-1}(\bfk^{-2} \mF(\bg(\bpsi^i_{f}))), ~\forall ~i=1,\ldots,n_f, \label{eq:poisson_for_basis}
\end{equation}
where $\bfk^{-2}$ is a vector containing the inverse of the squares of wavenumbers. In this approach, the potential reduced basis function  $\bpsi^i_\phi$ is explicitly related to the distribution reduced basis function $\bpsi^i_f$ through the Poisson problem. However, unlike in the first approach, $\bPsi_{\phi}$ is not necessarily orthonormal. 

\begin{remark}
    If the potential reduced space is built using the second approach, %In the second approach, 
    one does not have to solve a reduced Poisson problem to obtain the generalized coordinate of the potential $\hbphi$ because $\hbphi = \hbff$ in this case. In addition, the dimension of the potential subspace is the same as the dimension of the distribution subspace, that is, $n_\phi = n_f$. 
\end{remark}

For the electric field, the prescribed offset vector $\bE_\tos$ and the reduced
basis matrix 
\begin{equation}
    \bPsi_E =\begin{bmatrix} \bpsi^1_E & \bpsi^2_E & \cdots &
    \bpsi^{n_\phi}_{E} \end{bmatrix}
\end{equation}
are constructed with the 
relation $E = -\nabla \phi$, that is,
\begin{equation}
    \bE_\tos = \mF^{-1}(-{\bm i}\bfk\mF(
    \bphi_\tos
    )),\quad
    \bpsi^i_E = \mF^{-1}(-{\bm i}\bfk \mF(
    \bpsi^i_\phi
    )), \quad \forall~i=0,\ldots,n_\phi.
\end{equation}
with the gradient computed using FFT.

\subsection{Projection-based ROM}
\label{ssec:PROM}

With the solution bases introduced in the previous section %above reduced approximations 
and Galerkin projection, the reduced system of (\ref{eq:ode}) is written as\footnote{We note that because the reduced basis functions are orthonormal, the reduced mass
matrix $\bPsi_f^T \bPsi_f$ is an identity matrix.}:
\begin{align}
    \begin{cases}
   % \bPhi_f^T \bPhi_f 
    \displaystyle\frac{d \hbff}{dt} & = \bPsi_f^T \bN(\bff_\tos + \bPsi_f \hbff) \\[8pt]
    \wbA \hbphi & = \wbg,
    \end{cases}
    \label{eq:1D--1V-vlasov-rom}
\end{align}
where the reduced Poisson operator 
$\wbA\in \real^{n_\phi \times n_\phi}$ and the reduced right-handed side vector $\wbg \in \real^{n_\phi} $ are defined as: 
\begin{align}
    \wbA_{ij} & \equiv 
    (\bpsi^i_{\phi})^T \mF^{-1}(-\bfk^{2} \mF(\bpsi^j_{\phi})),\quad\text{and}\quad\wbg_i \equiv (\bpsi^i_{\phi})^T \left[\mF^{-1}(\bfk^{2} \mF(\bphi_{\tos})) + \bg({\bff_\tos} + \bPsi_f \hbff) \right] ,
\end{align} with $\bg(\bff)$ being the discretization of the right-hand side in (\ref{eq:1D--1V-poisson}). 
We note that FFT is used to compute the second
order derivative of the potential reduced basis functions $\bpsi^j_{\phi}$ for
constructing $\wbA$ and $\bfk^{2}$ is a vector containing the squares of
wavenumbers. In what follows, we refer (\ref{eq:1D--1V-vlasov-rom}) to the projection-based ROM.

\subsection{%Efficient treatment of nonlinear terms
Tensorial evaluation of nonlinear terms
}
\label{ssec:tensorial}

Although the projection-based ROM (\ref{eq:1D--1V-vlasov-rom}) is a system of ODEs with $n_f \ll N_f$ unknowns, the cost for evaluating the right-handed side of (\ref{eq:1D--1V-vlasov-rom}) scales with the FOM size, $N_f$, because it requires reconstructing the distribution field by multiplying the generalized coordinates $\hbff$ with the basis matrix $\bPsi_f$, computing $\bN(\bff_\tos+\bPsi_f \hbff)$ vector and projecting onto the reduced space $\mathcal{S}_f$. Therefore, we cannot
expect any speed-up with the projection-based ROM without special treatment of the nonlinear term.

{
To overcome this issue, we introduce a
tensorial approach to efficiently update the nonlinear term, which utilizes the tensor-product representation of the nonlinear function, $\bN$,
\begin{align}
    \bN(\bff) 
    & = \biggl( 
    [{\bDcx} \otimes {\bm I_v}] \left[{\bm I_x }\otimes \tdiag({\bm v}) \right] + [{\bm I_x} \otimes {\bDcv}]\left[ \tdiag(\bE) \otimes {\bm I_v} \right]\biggr) \bff, \label{eq:tensor_product} 
\end{align}
where $\bDcx$ and $\bDcv$ are the 1D conservative finite-difference operators correspond to (\ref{eq:cfd})--(\ref{eq:upwindflux}) in $x$ and $v$, respectively. $\tdiag(\bv)$ and $\tdiag(\bE)$ are the diagonal matrices with $\bv$ and $\bE$ values in the diagonal. 

With (\ref{eq:tensor_product}) and approximating $\bE$ with $\tbE$ (\ref{eq:reduced_app}), the nonlinear term in (\ref{eq:1D--1V-vlasov-rom}) can be approximated as
\begin{align}
    \bPsi_f^T \bN(\bff_\tos + \bPsi_f \hbff)
    & = \bPsi_f^T \left(
    [{\bDcx} \otimes \bIv] \left[\bIx \otimes \tdiag({\bm v}) \right] 
    + [\bIx \otimes {\bDcv}][ \tdiag(\bE_\tos + \sum_{i=1}^{n_\phi} \hbphi^i \bpsi^i_E) \otimes \bIv]
    \right) 
    \left(\bff_\tos + \bPsi_f \hbff \right) \nonumber \\ 
    & \approx \bPsi_f^T \biggl( 
   [{\bDcx} \otimes \bIv] \left[\bIx \otimes \tdiag({\bm v}) \right] + [\bIx \otimes {\bDcv}] \left[ \tdiag(\bE_\tos) \otimes \bIv \right] 
   \biggr) 
    \left( \bff_\tos + \bPsi_f \hbff \right) \nonumber \\ 
   & + \bPsi_f^T \biggl( [\bIx  \otimes {\bDcv}] \sum^{n_\phi}_{i=1} \hbphi^i  \left[ \tdiag(\bpsi^i_E) \otimes \bIv \right] 
   \biggr) 
    \left( \bff_\tos + \bPsi_f \hbff \right) \nonumber \\ 
     & {\approx}~ \wbG_1 \hbff + \bmG_2(\hbphi)\hbff + \wbg_0 + \wbG_0\hbphi + \wbg_1 + \wbG_2\hbff.
     \label{eq:tensorial_approx}
\end{align}
The vectors $\wbg_0 \in \real^{n_f} $, $\wbg_1 \in \real^{n_f}$, the matrices
$\wbG_0\in \real^{n_f \times n_\phi}$, $\wbG_1 \in \real^{n_f \times n_f}$, $\wbG_2 \in \real^{n_f \times n_f}$ and the tensor $\bmG_2 \in \real^{n_f \times n_f \times n_\phi}$ can be precomputed in the
offline and are defined as:
\begin{align}
    \wbG_{1,ij} & \equiv (\bpsi^i_{f})^T  \bigl[{\bDcx} \otimes \bIv \bigr] \bigl[{\bIx} \otimes \tdiag({\bv}) \bigr] \bpsi^j_{f} \\ 
    \bmG_{2,ijk} & \equiv (\bpsi^i_{f})^T \bigl[\bIx \otimes \bDcv \bigr] \bigl[ \tdiag(\bpsi^j_{E}) \otimes {\bIv} \bigr]  \bpsi^k_{f} \\ 
    {\wbg_{0,i}} & \equiv (\bpsi^i_{f})^T  \bigl[{\bDcx} \otimes \bIv \bigr] \bigl[{\bIx} \otimes \tdiag({\bv}) \bigr]  {\bff_{\tos}} \\ 
    \wbG_{0,ij} & \equiv (\bpsi^i_{f})^T \bigl[ \bIx \otimes {\bDv} \bigr]
    \bigl[ \tdiag(\bpsi^j_{E}) \otimes {\bIv} \bigr]  {\bff_{\tos}} \\
    {\wbg_{1,i}} & \equiv (\bpsi^i_{f})^T    
    \bigl[ \bIx \otimes \bDv \bigr] \bigl[ \tdiag(\bE_\tos) \otimes \bIv \bigr] {\bff_{\tos}} \\
    \wbG_{2,ij} & \equiv (\bpsi^i_{f})^T  
    \bigl[ \bIx \otimes \bDv \bigr] \bigl[ \tdiag(\bE_\tos) \otimes \bIv \bigr]\bpsi^j_{f}.
\end{align}
We note that, in addition to the reduced basis approximation to the electric field, there is another approximation in (\ref{eq:tensorial_approx}) due to the nonlinear operators $\bDcx$ and $\bDcv$. These 1D conservative finite-difference operators are nonlinear because the WENO/Upwind scheme is used to evaluate the left- and right-biased approximations of the flux at the interface. 
}
The tensorial ROM is defined to be the following system of ODEs,
\begin{align}
    \begin{cases}
    \displaystyle\frac{d \hbff}{dt} & =
    \wbG_1 \hbff + \bmG_2(\hbphi)\hbff + \wbg_0 + \wbG_0\hbphi + \wbg_1 + \wbG_2\hbff \\[8pt]
    \wbA \hbphi & = \wbg
    \end{cases}
    \label{eq:tensorial-rom}.
\end{align}
\begin{remark}
    The cost for solving the tensorial ROM (\ref{eq:tensorial-rom}) is independent of the FOM degrees of freedoms $N_f$ and scales as the cost of tensor contraction $\mathcal{O}(n^2_fn_\phi)$.
\end{remark}

\section{Indicator-Based Decomposition of Solution Manifold}
\label{sec:temporal_ROM}

Although the cost for solving the tensorial ROM introduced in
Section~\ref{sec:ROMs} is independent of the FOM degrees of freedoms, the
advection-dominated nature of the solution of the simulation of Vlasov\dash Poisson
equations requires large reduced space dimensions in order for ROMs to be
accurate. For example, in the two-stream
instability problem, there are three stages of the solutions, that is, short transient, growth, and statistically stationary. As a result, the linear dependence among the snapshots is weak, and therefore there is no intrinsic low-dimensional subspace that can approximate the solution manifold comprised of all the solutions over the temporal domain. In addition, the cost of solving the tensorial ROM scales like
$\mathcal{O}(n^2_fn_\phi)$ due to the tensor contraction cost per time step. Hence, it is impossible for the reduced-order model approach to achieve any meaningful speed-up and good accuracy.

To overcome these difficulties, we consider employing multiple
reduced-order models in time. The idea of the methodology is to construct local ROMs
in the parameter-time domain using a suitable indicator %physical time $t$ as
for clustering and
classification~\cite{copeland2022reduced,cheung2023local}. In the offline phase, we construct each of these ROMs from a small subset of the snapshot samples to ensure low dimension. In the online phase, each of these ROMs is used in a certain subset of the parameter-time domain where they are supposed to provide a good approximation. The idea is to decompose the solution manifold into submanifolds where the Kolmogorov $n$-width decays fast with respect to the subspace dimension. This enables us to build accurate multiple low-dimensional subspaces.
We note that the methodology is not limited to the 1D1V Vlasov\dash Poisson equations and has been applied to the Euler equations \cite{cheung2023local} and the Navier\dash Stokes equations \cite{parish2021windowed}. In addition, the framework is not restricted to the tensorial ROM (\ref{eq:tensorial-rom}) and can be applied to other ROMs, e.g., DMD-based ROM.

We consider two indicators for decomposition, namely, the physical time and the electric field energy. We note that evaluating both indicators does not lower the speed-up of the overall ROM simulation in the online stage because the costs of the indicator evaluation do not involve operations that scale with FOM.

\subsection{Time-windowing: Decomposition by physical time}

A very natural choice of indicator is the physical time, that is,
$\indicator(\bff,t,\bp) = t$, where $\indicator:\real^{N_f} \times \real^+
\times \mathcal{D} \rightarrow \real$ is a function that maps the triplet
$(\bff,t,\bp)$ to a real value. In the following, we describe the framework of
the physical-time decomposition of the solution manifold. 

To begin with, the range of the indicator
$[\indicator_{\textup{min}},\indicator_{\textup{max}})$, is partitioned into
$N^{\textup{off}}_w$ subintervals, namely,
\begin{equation}
    \indicator_{\textup{min}} = \indicator_0 < \indicator_1 < \cdots < \indicator_{N^\textup{off}_w -1} < \indicator_{N^\textup{off}_w} = \indicator_{\textup{max}}.
\end{equation}
This partitioning can be either prescribed or determined by the snapshot data collected in the offline phase. With the partitioning of the indicator range, instead of directly assembling all the snapshot samples into a single huge snapshot matrix (\ref{eq:full_snapshot_matrix}), the FOM states are first classified into groups. Let $m \in \mathbb{N}(N^\textup{off}_w)$ be a group of indices. We denote the subset of paired indices of time and the offline parameter whose corresponding snapshot belongs to the $m$-th group as

\begin{equation}
   \mathcal{G}_m = \{(n,k) \in \mathbb{Z} \times \mathbb{N}(n_{\bp}): 0\le n \le N_t(\bp) ~\text{and}~ \indicator(\bff_n(\bp_k), t_n, \bp_k) \in [\indicator_{m-1}, ~\indicator_m) \},\label{eq:group}
\end{equation}
where $N_t(\bp)$ is the number of time steps used in the FOM simulation for parameter $\bp$, and the final time $t_f = t_{N_t(\bp)}$ is considered to be uniform over the parameter domain $\mD$. 
Then the distribution snapshot matrix $\bF_m$ in the $m$-th group is formed by assembling the corresponding snapshots, that is, 
$\bF_{m} \equiv \left[ \bff_n(\bp_k) \right]_{(n,k) \in \mathcal{G}_m}$ and POD is used to construct the distribution solution basis matrix $\bPsi^m_{f}$ as described in Section \ref{sec:solution_basis_construction}. We note that the $m$-th group $\mG_m$ is also used to form the potential snapshot matrix $\bF_{\phi,m}$, that is $\bF_{\phi,m} \equiv \left[ \bphi_n(\bp_k) \right]_{(n,k) \in \mathcal{G}_m}$.

For a generic problem parameter $\bp \in \mathcal{D}$, the computation in the online phase is performed using different reduced bases in $N_w$ subintervals of the temporal domain $[0,~t_f]$, that is, 
\begin{equation}
    0 = T_0(\bp)< T_1(\bp)< \cdots < T_{N_w(\bp)-1}(\bp) < T_{N_w(\bp)}(\bp) = t_f(\bp),
\end{equation}
with $\{T_j(\bp)\}^{N_w(\bp)-1}_{j=1}$ being the partition of the temporal domain determined by the indicator $\indicator$. 
In general, the temporal domain partition is parameter-dependent, that is, the number of subintervals $N_w(\bp)$ and the end-point of each subinterval $\{T_j(\bp)\}^{N_w(\bp)}_{j=1}$ depend on the problem parameter $\bp$, and they are assigned by the indicator $\indicator$ iteratively in the ROM online stage. 
We note that in the case of the physical time indicator, the temporal domain partition is parameter-independent, and 
for any problem parameter $\bp \in \mathcal{D}$, we have $T_j(\bp) =
\indicator_j$ for all $1 \le j \le N_w(\bp)$. In addition, $N_w(\bp) = N^{\textup{off}}_w$ because 
the final time is uniform over the parameter domain $\mathcal{D}$. We note that this choice reduces to the time-windowing (TW) reduced-order model approach in \cite{copeland2022reduced}.

Therefore, for $t \in \mathcal{T}_j\equiv [T_{j-1}(\bp) =\indicator_{j-1},~T_j(\bp)=\indicator_{j}]$, we employ the reduced bases constructed from the snapshot group $\mathcal{G}_m$. More precisely, we use the solution representation
\begin{align}
   \tbff(t;\bp) = \bPsi^m_{f} \hbff^j(t; \bp), \quad 
   \tbphi(t;\bp) = \bPsi^m_{\phi} \hbphi^j(t; \bp), 
\end{align}
where $\bPsi^m_f$ and $\bPsi^m_{\phi}$ are 
distribution and potential solution basis matrices constructed using POD from
$\bF_{m}$ and $\bF_{\phi,m}$, respectively. 
$\hbff^j:\mathcal{T}_j \times \mathcal{D} \rightarrow \mathbb{R}^{n_f}$ and $\hbphi^j:\mathcal{T}_j \times \mathcal{D} \rightarrow \mathbb{R}^{n_\phi}$ are the time-dependent generalized coordinates for the distribution and potential fields in the time interval $\mathcal{T}_j$, respectively. 
Following the same derivation in Section~\ref{ssec:tensorial}, we obtain the time-windowing ROM (TW-ROM), where
\begin{equation}
    \begin{cases}
    \displaystyle\frac{d \hbff^j}{dt} &=
    \wbG^m_1 \hbff^j + \bmG^m_2(\hbphi^j)\hbff^j + \wbg^m_0 + \wbG^m_0\hbphi^j + \wbg^m_1 + \wbG^m_2\hbff^j\\[8pt]
     \wbA^m \hbphi^j & = \wbg^m
    \end{cases}
    \label{eq:TW--ROM}
%   \label{eq:temporal_rom}.
\end{equation}
is a tensorial ROM for the time interval $\mT_j$ and the initial condition is given by projecting the solution at $t= T_{j-1}(\bp)$ in the time interval $\mT_{j-1}$ onto the reduced space for $\mT_{j}$, that is, 
\begin{equation}
    \hbff^j(T_{j-1}(\bp); \bp) = (\bPsi^m_f)^T \tbff(T_{j-1}(\bp);\bp).
\end{equation}

\subsection{Energy-windowing: Decomposition by electric field energy}
\label{ssec:EW--ROM}

\begin{figure}[!ht]
    \centering
    \includegraphics[width=0.9\textwidth]{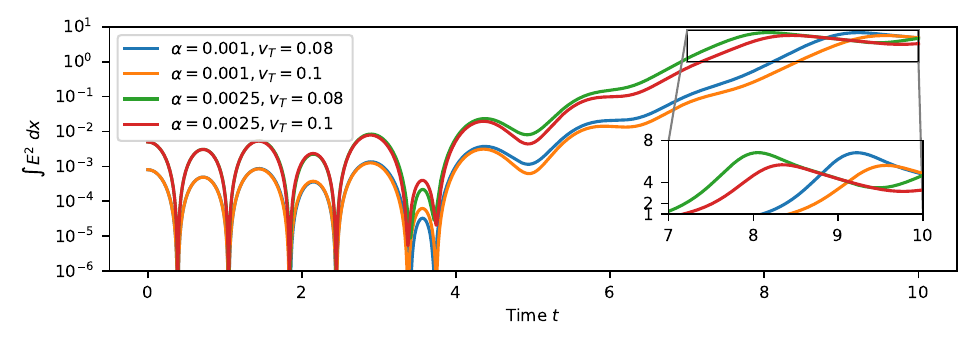}
    \caption{Two-stream instability. The behavior of the FOM electric field energy $\int E^2~dx$ in time for parameters $\bp = (\alpha, \vT) =(0.0025,0.08)$, $(0.0025,0.1)$, $(0.001,0.08)$, and $(0.001,0.1)$.} %\pht{Move the plot up to Section 4.}}
    \label{fig:fom_maxE}
\end{figure}

The evolution of the electric field energy $\int E^2~dx$ for two-stream instability (Section~\ref{ssec:two_stream}) is shown in Fig.~\ref{fig:fom_maxE} 
for four parameters $\bp \in \mD$, namely $\bp = (0.0025,0.08),~(0.0025,0.1),~(0.001,0.08),~\text{and}~(0.001,0.1)$. 
The results indicate that there are three stages of the solutions, namely, short transient, growth, and statistically stationary stages.
The short transient stage and the beginning time of the growth stage are parameter-independent. However, the end time of the growth stage is parameter-dependent. With the physical time indicator, the temporal domain partition is parameter-independent, which could lead to having  
snapshot samples that are very different but in the same group \cite{cheung2023local}, and prevent the linear subspace method from being effective and limit the speed-up.

To improve the classification of snapshots in the growth and statistically stationary stages for different parameters, and further reduce the size of the solution basis, we propose to use the electric field energy, which characterizes the solutions of the Vlasov\dash Poisson equations, as an indicator to decompose the solution manifolds in those two stages. The electric field energy indicator is defined as 
\begin{equation}
    \indicator(\bff,t,\bp) =\int_{\Omega_x} \bE^2~dx,
\end{equation}
where $\bE$ is the electric field solution associated with the distribution solution $\bff$ at time $t$ and parameter $\bp$. 

With the electric field energy indicator, the temporal domain is partitioned into $N^{\text{off}}_w =5$ subintervals. In the offline stage, the FOM states are classified into our five groups defined below:
\begin{align}
   \mG_1 & = \{(n,k) \in \mathbb{Z} \times \mathbb{N}(n_{\bp}): 0\le n \le N_t(\bp) ~\text{and}~ t_n \in [0, ~4) \} \\ 
   \mG_2 & = \{(n,k) \in \mathbb{Z} \times \mathbb{N}(n_{\bp}): 0\le n \le N_t(\bp), ~ t_n \ge 4~\text{and}~ \indicator(\bff_n(\bp_k),t_n,\bp_k) \le 1\} \\ 
   \mG_3 & = \{(n,k) \in \mathbb{Z} \times \mathbb{N}(n_{\bp}): 0\le n \le N_t(\bp), ~ \indicator(\bff_n(\bp_k),t_n,\bp_k) > 1~\text{and}~\frac{\partial \indicator}{\partial t}(\bff_n(\bp_k),t_n,\bp_k) \ge 0\} \\
   \mG_4 & = \{(n,k) \in \mathbb{Z} \times \mathbb{N}(n_{\bp}): 0\le n \le N_t(\bp), ~ (n,k) \notin \mG_1 \cup \mG_2 \cup \mG_3 ~\text{and}~\frac{\partial \indicator}{\partial t}(\bff_n(\bp_k),t_n,\bp_k) \le 0\} \\
   \mG_5 & = \{(n,k) \in \mathbb{Z} \times \mathbb{N}(n_{\bp}): 0\le n \le N_t(\bp), ~\text{and}~(n,k) \notin \mG_1 \cup \mG_2 \cup \mG_3 \cup \mG_4 \}.
\end{align}
Then the snapshot matrices $\bF_m$ and $\bF_{\phi,m}$ in the $m$-th group are formed by assembling the corresponding snapshots
and POD is used to construct the reduced basis matrices $\bPsi^m_{f}$ and $\bPsi^m_{\phi}$.

For a generic problem parameter $\bp \in \mathcal{D}$, the computation in the online phase is performed using different reduced bases in the five subintervals of the temporal domain $[0,~t_f]$, that is, 
\begin{equation}
    0 = T_0(\bp)< T_1(\bp) = 4 < T_2(\bp) < T_{3}(\bp) < T_{4}(\bp) < T_{5}(\bp) = t_f(\bp).
\end{equation}
We note that the temporal domain partition is parameter-dependent because the
end-points $T_2(\bp)$, $T_3(\bp)$ and $T_4(\bp)$ are determined by the electric field energy
indicator and are assigned by the indicator iteratively in the ROM online stage.

Therefore, for $t \in \mathcal{T}_j\equiv [T_{j-1}(\bp),~T_j(\bp)]$, we employ the reduced bases constructed from the snapshot group $\mathcal{G}_m$ and follow the same derivation in Section~\ref{ssec:tensorial}, we obtain the energy-windowing ROM (EW-ROM), where 
\begin{equation}
    \begin{cases}
        \displaystyle \frac{d \hbff^j}{dt} & =
        \wbG^m_1 \hbff^j + \bmG^m_2(\hbphi^j)\hbff^j + \wbg^m_0 + \wbG^m_0\hbphi^j + \wbg^m_1 + \wbG^m_2\hbff^j
        \\[8pt]
        \wbA^m \hbphi^j & = \wbg^m
    \end{cases}
    \label{eq:EW--ROM}
\end{equation}
is a tensorial ROM for the time interval $\mT_j$ and the initial condition is given by projecting the solution at $t= T_{j-1}(\bp)$ in the time interval $\mT_{j-1}$ onto the reduced space for $\mT_{j}$, that is, 
\begin{equation}
    \hbff^j(T_{j-1}(\bp); \bp) = (\bPsi^m_f)^T \tbff(T_{j-1}(\bp);\bp).
\end{equation}

\begin{figure}[!ht]
    \centering
    \includegraphics[width=0.9\textwidth]{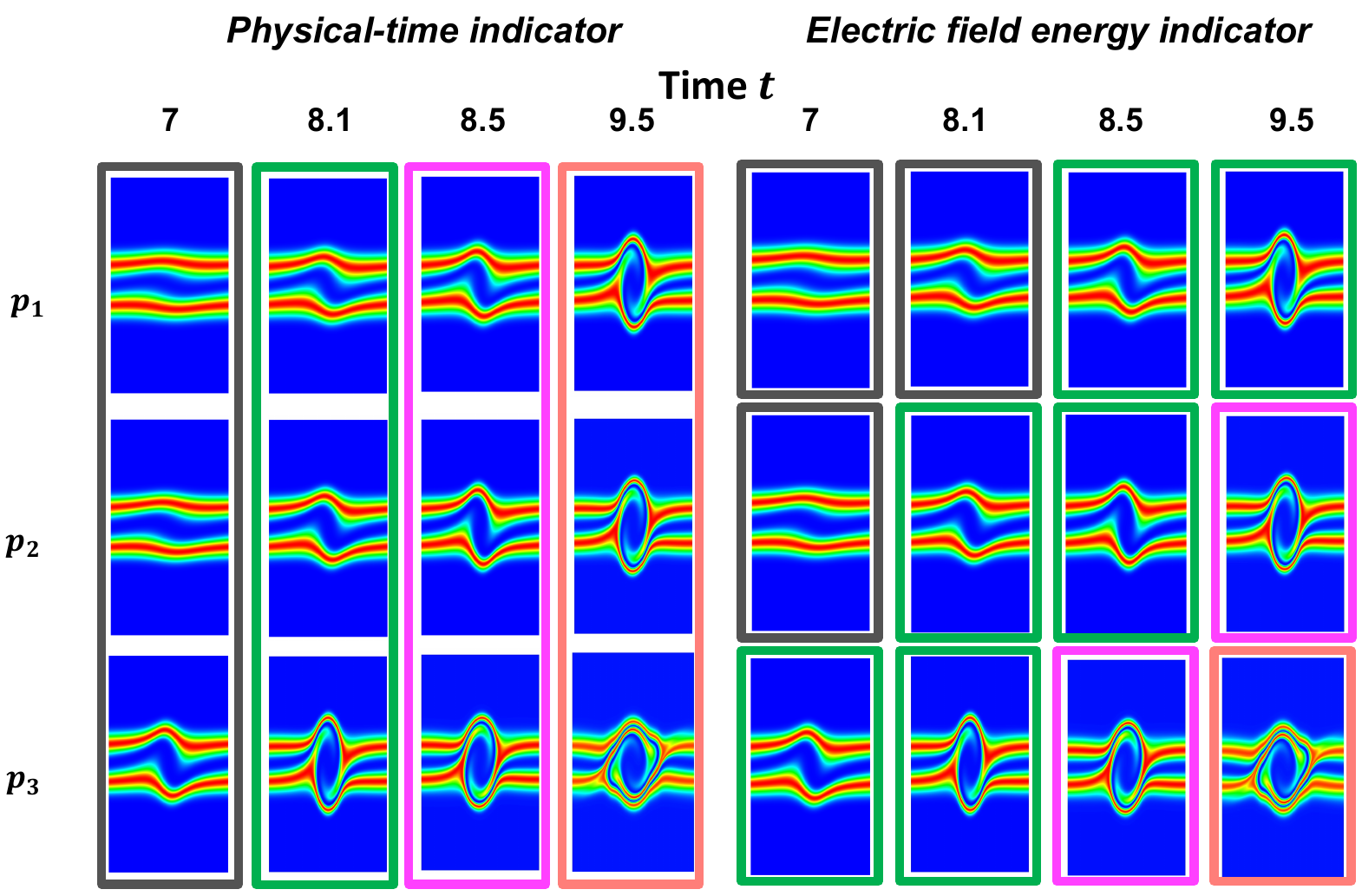}
    \caption{An illustrative example to explain the mechanism of decomposition
    of solution manifold using the physical time and the electric field energy as
    indicator. The samples in each group are surrounded by a box with the same
    color, namely dark gray, green, green, pink, and orange. Groups classified using the electric field energy  
    have strong linear dependence.}
    \label{fig:}
\end{figure}

\section{Numerical Results}
\label{sec:numerical}

In this section, we present numerical results to test the performance of our proposed method applied to three benchmark problems of the 1D1V Vlasov\dash Poisson equations. %the LS-ROM with tensorial and time-windowing approaches. 
All the simulations in this section use Intel Sapphire CPUs on the machine Dane in the Livermore Computing Center \footnote{High-performance computing at LLNL, \url{https://hpc.llnl.gov/hardware/compute-platforms/dane}}, with $256$ GB memory and peak TFLOPS (CPUs) of $10,723$. 

To evaluate the ROM performance, we consider the relative error between the FOM and the ROM distribution field at the final time $t_f$, which is defined as:
\begin{align}
   \epsilon_{f,{\tilde{t}_f}} \coloneqq \frac{\|\bff - \tbff \|_2}{\|\bff\|_2}.
\end{align}
The speed-up of each ROM simulation is measured by dividing the wall-clock time for the FOM time loop by the wall-clock time for the corresponding ROM time loop.

\subsection{2D (1D1V) Vlasov equation - Prescribed $E$-field} 
We first consider a model problem in which the effects of the self-consistent
electric field $E$ are negligible compared to the ones caused by an external
electric field $E_{\text{prescribed}}(x) = \alpha \cos(x)$, which is 
independent of time and periodic in the $x$-direction. The initial solution is:
\begin{align}
   f(x,v) = \left[ 1 + \alpha \cos\left(x\right) \right]
   \exp\left[-\frac{1}{2}\left(\frac{v}{\vT}\right)^2\right],
\end{align}
where the perturbation amplitude $\alpha$ and the thermal velocity $\vT$ are fixed to $0.1$ and $1$, respectively. The computational domain is $(x,v) \in \Omega \coloneqq [0,2\pi]\times
[-6v_T, 6v_T]$ with periodic boundary condition in the $x$-direction
and homogeneous Dirichlet in the $v$ direction. 
The FOM has $N_f = 128 \times 128 = 16384$ degrees of freedoms and is simulated in the time interval $\Tinterval = (0, t_f]$, with $t_f= 130$ and $\Delta t=0.005$. 

We consider the reproduction problem,  where the ROM is tested at the same
parameter $\bp$ and time interval $\Tinterval$ as of the FOM. The reduced basis
matrix $\bPsi_f$ is built using $K=2600$ snapshots collected in the time interval $\Tinterval$ and
further used to construct the projection-based ROM (\ref{eq:1D--1V-vlasov-rom}). 
\begin{remark}
    For this model problem, the tensorial ROM (\ref{eq:tensorial-rom}) is not required because 
    the prescribed electric field $E_{\text{prescribed}}(x)=\alpha\cos\left(x\right)$ is independent of the solution; therefore, the equation is linear.
\end{remark}

\begin{figure}[!ht]
    \centering
    \includegraphics[width=0.9\textwidth]{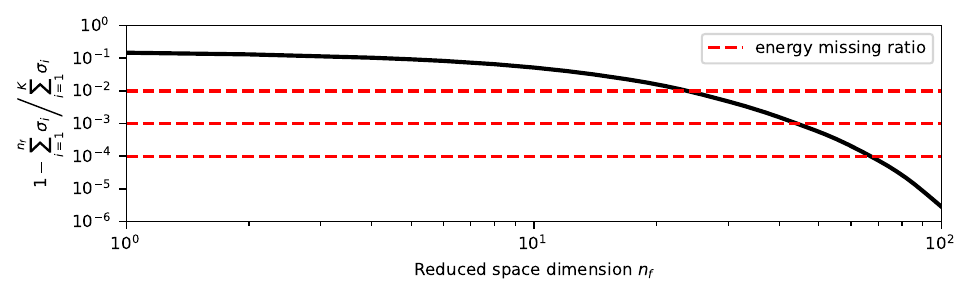}
    \caption{Prescribed electric field case. %with perturbation amplitude $\alpha=0.1$ and final time $t_f=130$: 
    %(a) Singular values of the snapshots matrix $\bF$. (b) 
    The behavior of $1-\sum^{n_f}_{i=1} \sigma_i /\sum^K_{i=1} \sigma_i$ as a function of the reduced space dimension $n_f$.}
    \label{fig:preE_singular_T130}
\end{figure}

The reducibility of the problem is studied by computing singular values of the snapshots matrix $\bF \in \mathbb{R}^{N_f \times K}$. 
\begin{figure}[!ht]
    \centering
    \includegraphics[width=0.9\textwidth]{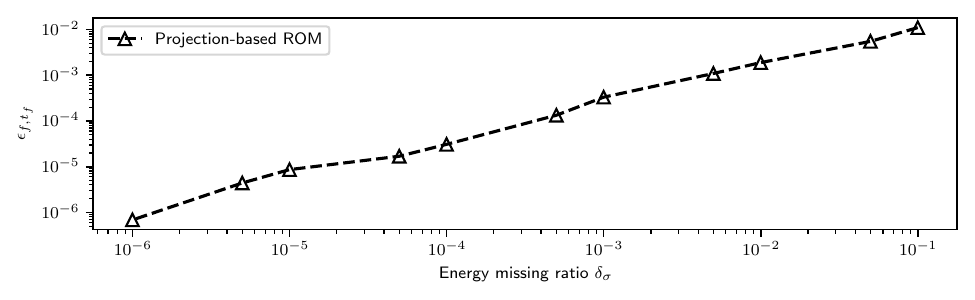}
    \caption{Prescribed electric field case. The behavior of the relative error in the distribution field at the final time $t_f$, $\epf$ with respect to energy missing ratio $\delta_\sigma$.}%with perturbation amplitude  $\alpha=0.1$ and final time $t_f=130$.} 
    \label{fig:pE_reprod_error}
\end{figure}
In Fig.~\ref{fig:preE_singular_T130}, 
the behavior of $1 -{\sum^{n_f}_{i=1}\sigma_i}/{\sum^{K}_{i=1}\sigma_i}$ with respect to the reduced space dimension $n_f$ is shown 
with three energy missing ratio values, that is, $\delta_\sigma = \tento{-2},~\tento{-3}$ and $\tento{-4}$,  represented by the red dashed line.
Using (\ref{eq:energy_criteria}), we find that $n_f=24$ reduced basis functions are able to capture $99\%$ of the energy in the snapshots (i.e., $\delta_\sigma=\tento{-2}$), indicating that the global linear subspace can effectively represent the solutions in the time interval $\Tinterval=(0,130]$.

\begin{figure}[!ht]
     \centering
     \resizebox{0.6\columnwidth}{!}{%
    \includegraphics[width=0.49\textwidth]{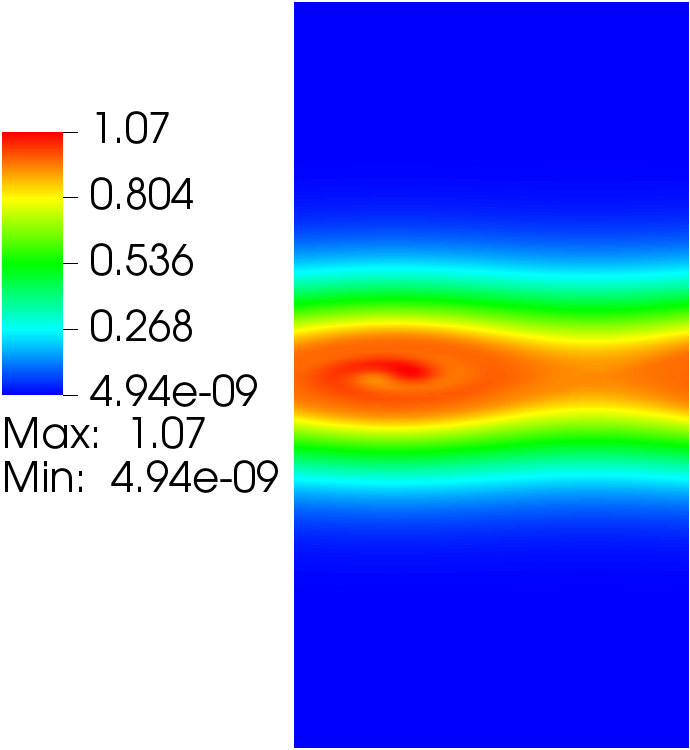}
    \hspace{1cm}
    \includegraphics[width=0.49\textwidth]{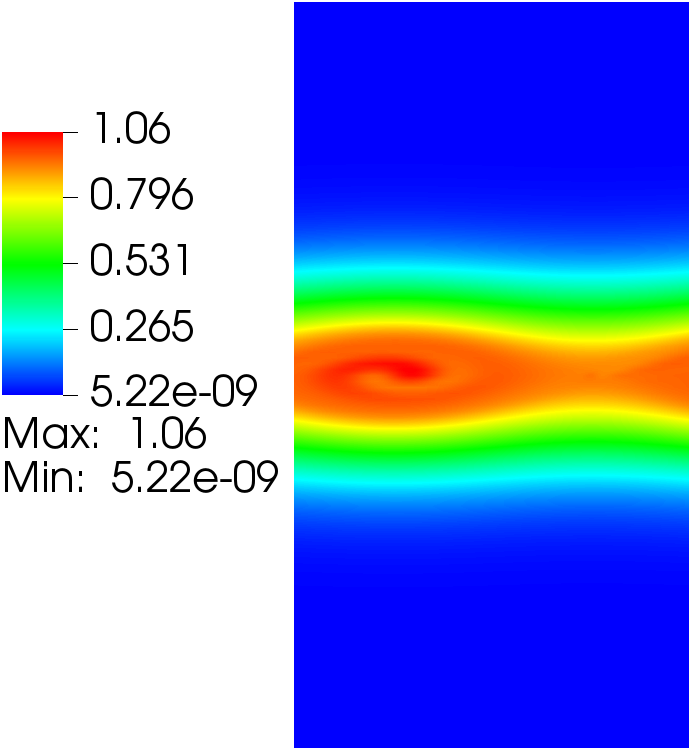}
    }
    \caption{Prescribed electric field case. %From left to right: Initial condition, 
    The FOM and ROM solutions at time $t=130$ with energy missing ratio $\delta_\sigma = \tento{-2}$.}
    \label{fig:preE_init_final}
\end{figure}

Fig.~\ref{fig:pE_reprod_error} shows 
the relative error in the distribution field at the final time $t_f$, $\epf$, of the 
projection-based ROM with various energy missing ratio $\delta_\sigma$.
The results show that $\epf$ decreases
as $\delta_\sigma$ decreases and is able to reach to $\tento{-6}$ with $\delta_\sigma=\tento{-6}$. 
We note that already with
$\delta_\sigma= \tento{-2}$, the projection-based ROM is able to achieve an error of $0.2\%$ with a speed-up of $10$ and the reconstructed solution at the final time $t_f=130$ along with the FOM solution are shown in Fig.~\ref{fig:preE_init_final}. 

\subsection{2D (1D1V) Vlasov\dash Poisson equation - Landau damping}

We consider the Landau damping \cite{finn2023numerical} case, which is a model
problem used to study the damped propagation of small amplitude plasma waves.
Landau damping plays a significant role in plasma physics and can be applied to
study and control the stability of charged beams in particle accelerators.
The initial condition is
\begin{align} 
      f(x,v) = \frac{1}{\sqrt{2\pi \vT^2}}\exp\left[-\frac{1}{2}\left(\frac{v}{\vT}\right)^2\right] \left[ 1+\alpha \cos (k x)\right] = f_0(v) + \alpha \cos(kx) f_0(v),
\end{align}
where the perturbation wavenumber $k$ is fixed to $1$ and the perturbation amplitude $\alpha$ and 
the thermal velocity $\vT$ are the study parameters $\bp = (\alpha, \vT) \in \mD = [0.01, 0.03]\times[0.8,1]$. We consider an offset $f_\tos(v)$ for the ROM, which is defined as  
\begin{equation}
    f_\tos(v) = f_0(v)   \frac{1}{\sqrt{2\pi \vT^2}}\exp\left[-\frac{1}{2}\left(\frac{v}{\vT}\right)^2\right].
\end{equation}
Because of the choice of offset for the distribution field, the offset for the potential and electric fields, $\bphi_\tos$ and $\bE_\tos$ are zero vectors. 

\begin{figure}[!ht]
    \centering
     \resizebox{0.8\columnwidth}{!}{%

    \includegraphics[width=0.24\textwidth]{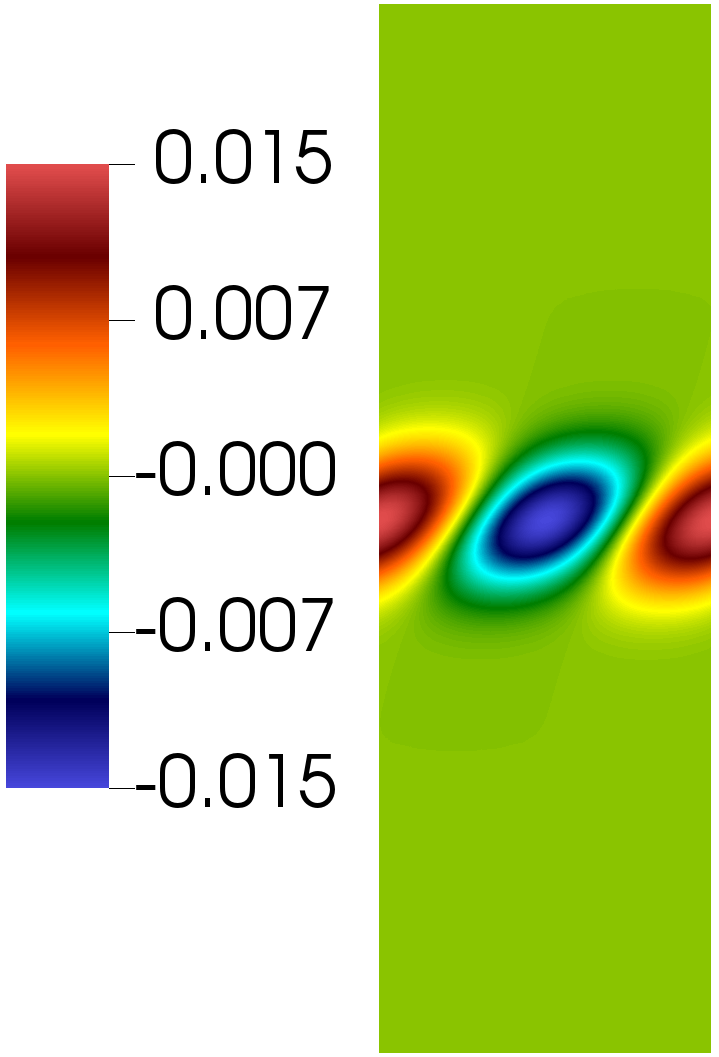}
    \includegraphics[width=0.24\textwidth]{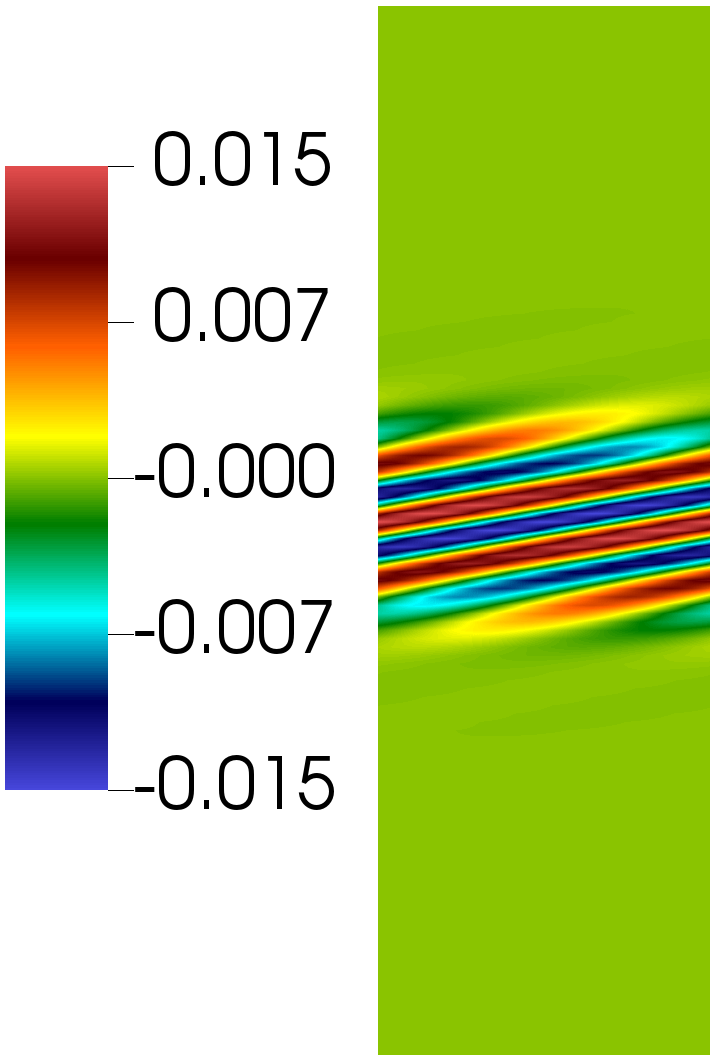}
    \includegraphics[width=0.24\textwidth]{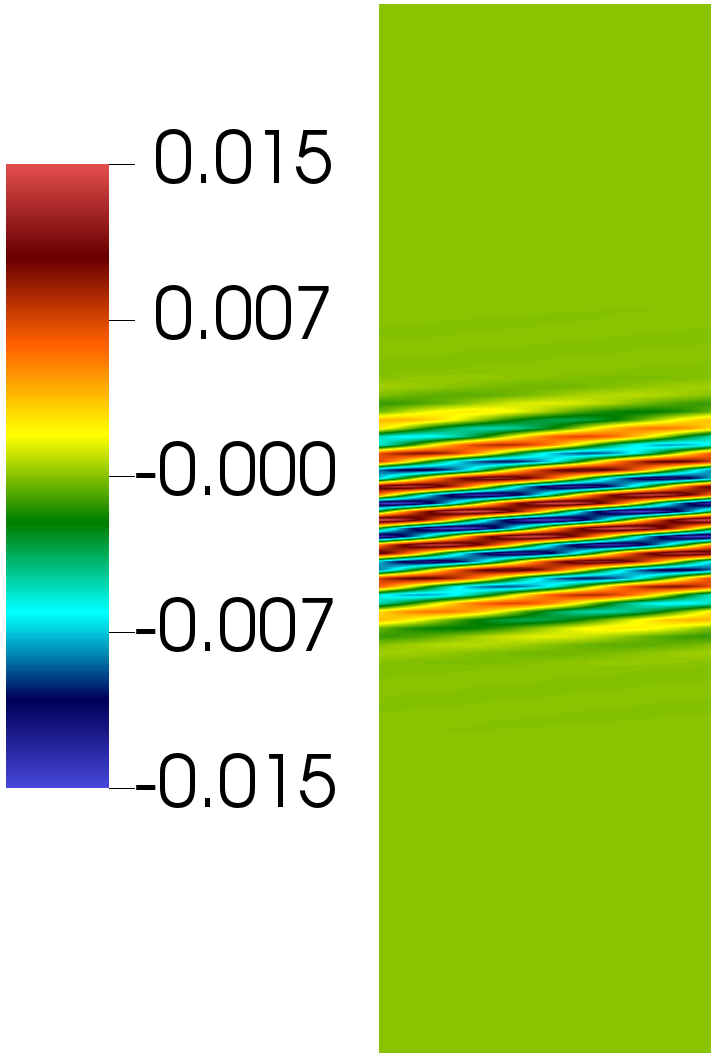}
    \includegraphics[width=0.24\textwidth]{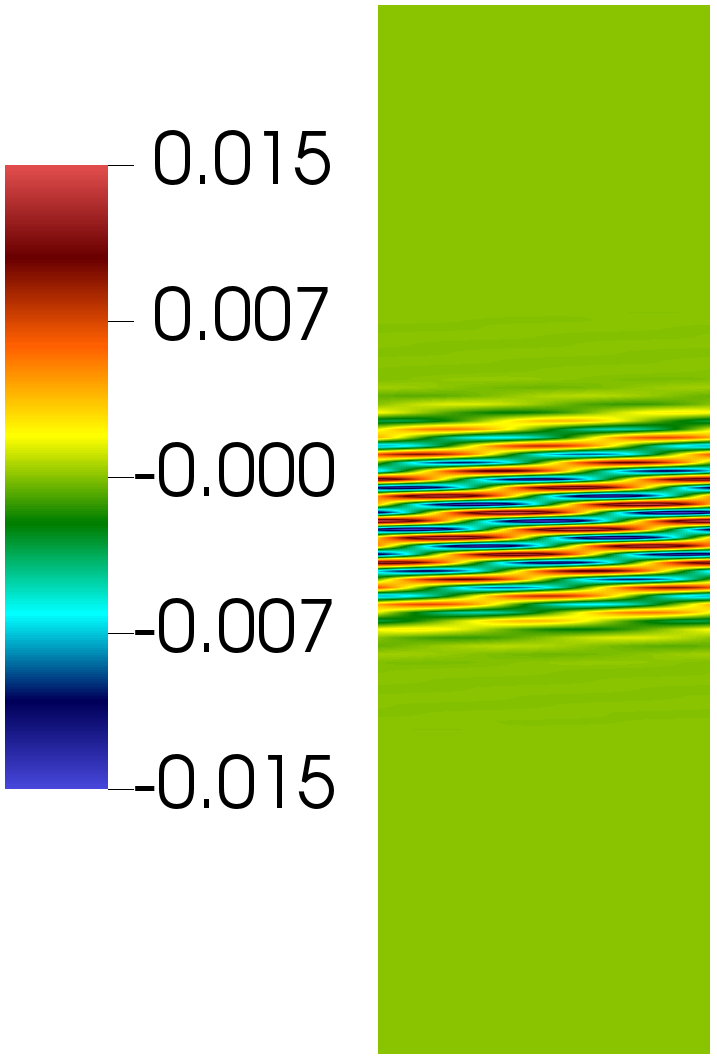}
    }
    \caption{Landau damping. Left to right: The FOM perturbation distribution fields at time $t=0$, $t=5$, $t=10$ and $t=15$.} 
    \label{fig:ld_fom_snapshot_centered}
\end{figure}

We consider the computational domain $(x,v) \in \Omega \coloneqq [0,2\pi/k]\times [-10, 10]$ with periodic boundaries in the $x$-direction and homogeneous Dirichlet boundaries in the $v$-direction. 
The FOM has $N_f = 128 \times 128 = 16384$ degrees of freedoms. A uniform time step $\Delta t = 0.005$ is used for the evolution of the particle's distribution over the time interval
$\Tinterval=(0, t_f]$ with $t_f = 15$. 
The centered FOM distribution fields $\bff - \bff_\tos$ with $\bp = (0.03, 0.8)$ 
at time $t=0,5,10$ and $15$ are shown in Fig.~\ref{fig:ld_fom_snapshot_centered}.

\subsubsection{Reproduction problem}

We consider reproduction problems at $\bp = (0.01, 1)$ and $\bp = (0.03, 0.8)$, which correspond to the simplest and the most challenging cases in the parameter space $\mD$. In each reproduction problem, the single ROM (\ref{eq:tensorial-rom}) is constructed by performing POD on $K=3000$ snapshots collected in the time interval $\Tinterval$. 
For the TW-ROM (\ref{eq:TW--ROM}), the range of the 
physical time indicator $[\indicator_{\text{min}} = 0,~
\indicator_{\text{max}}=15)$ is partitioned uniformly into $N_w$ subintervals. 
$K=3000$ snapshots are then classified into groups based on the partition of the indicator range. We remark that when $n_{p} = 1$, according to (\ref{eq:group}), the $j$-th group is  
\begin{align}
   \mG_j = \{(n,1) \in \mathbb{Z} \times \{1\}: 0\le n \le N_t ~\text{and}~ t_n \in [\indicator_{j-1}, ~\indicator_j)\}.
\end{align}
Then, for $1 \le j \le N_w$, the TW-ROM is constructed by performing POD on the subset of snapshots whose indices belong to the group $\mG_j$. In the online phase, the corresponding reduced order model is used, where the end-point $T_j(\bp)$ of the temporal subinterval $\mT_j(\bp)$ is defined as the time instance $\widetilde{t_n}(\bp)$ when $\indicator(\tbff_n(\bp_k), \widetilde{t_n}, \bp_k)$ first exceeds $\indicator_j$, at which we increment to the next subinterval $\mT_{j+1}(\bp)$.
\begin{figure}[!ht]
    \centering
    \begin{subfigure}{1\columnwidth}
    \centering
    \caption{$\bp = (0.01, 1)$}
    \includegraphics[width=0.9\textwidth]{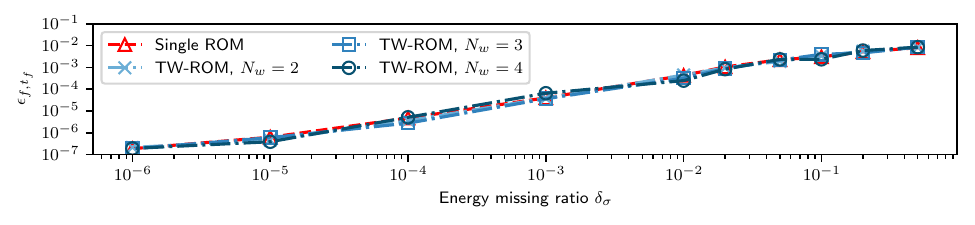}
    \end{subfigure}
    \begin{subfigure}{1\columnwidth}
    \centering
    \caption{$\bp = (0.03,0.8)$}
    \includegraphics[width=0.9\textwidth]{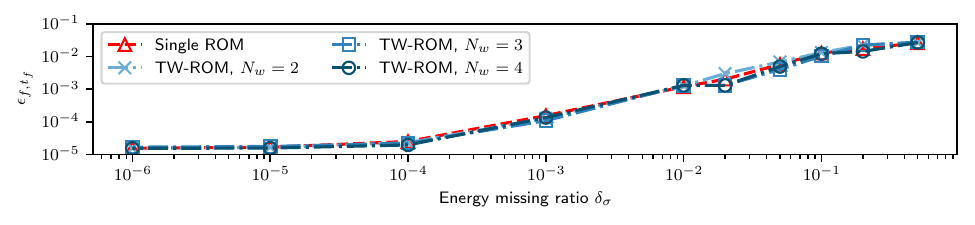}
    \end{subfigure}
    \caption{Reproduction problems of the Landau-damping problem with final time $t_f=15$ for parameters $\bp = (0.01,1)$ and $\bp=(0.03,0.8)$. The relative error in the distribution field at the final time $\epf$ of the single ROM and the TW-ROMs with respect to the energy missing ratio $\delta_\sigma$ are shown.} 
    \label{fig:ld_reprod_error_compare}
\end{figure}

In Fig.~\ref{fig:ld_reprod_error_compare},
the behavior of the relative error $\epf$ with respect to $\delta_\sigma$ for $\bp=(0.01,1)$ and $\bp=(0.03,0.8)$ are shown.
The results show that both $\epf$ of the single ROM and the TW-ROMs decrease as $\delta_\sigma$ decreases. Additionally, 
the error of the single ROM is comparable to that of the TW-ROMs for both $\bp$.
Furthermore, we find that increasing the number of subintervals $N_w$ in the TW-ROM does not further improve accuracy.

\begin{figure}[!ht]
    \centering
    \begin{subfigure}{0.8\columnwidth}
    \centering
    \caption{$\bp = (0.01, 1)$}
    \includegraphics[width=1\columnwidth]{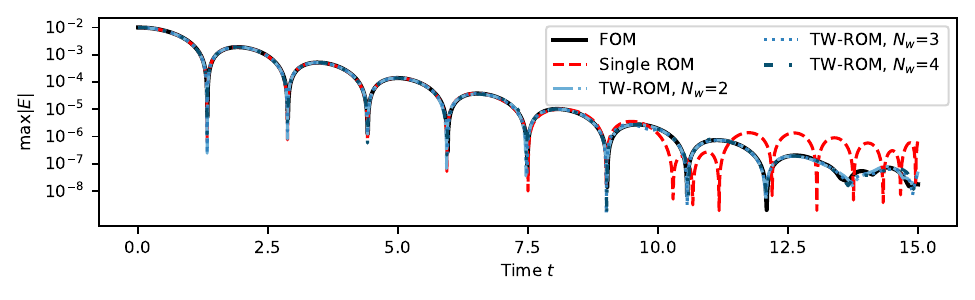}
    \end{subfigure}
    \begin{subfigure}{0.8\columnwidth}
    \centering
    \caption{$\bp = (0.03, 0.8)$}
    \includegraphics[width=1\columnwidth]{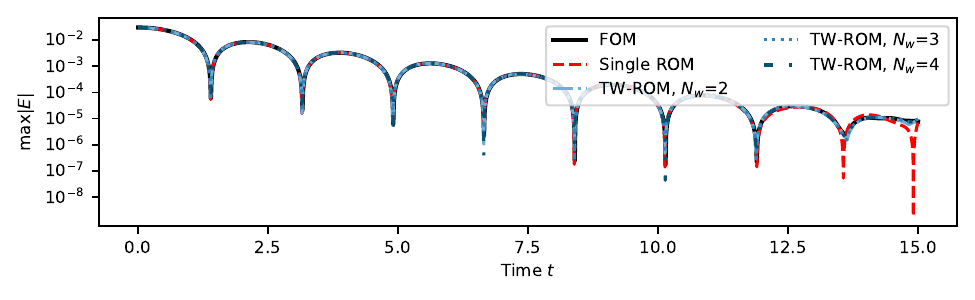}
    \end{subfigure}
    \caption{Reproduction problem of the Landau damping for parameters $\bp=(0.01, 1)$
    and $\bp=(0.03,0.8)$. The behavior of the $\max |E|$ in the time interval $\Tinterval$ for
    the single ROM
    (\ref{eq:tensorial-rom}) and the TW-ROM (\ref{eq:TW--ROM}) with the energy missing ratio $\delta_\sigma=\tento{-3}$.} %\pht{TODO: replot with upw5 results.}}
    \label{fig:ld_reprod1_maxe}
\end{figure}

We fix the energy missing ratio $\delta_\sigma$ to be $\tento{-3}$ for all the
ROMs and further examine the approximated $\max |E|$ in the time interval
$\Tinterval$ for $\bp=(0.01,1)$ and $\bp=(0.03,0.8)$ in
Fig.~\ref{fig:ld_reprod1_maxe}. 
We observe that the single ROM is able to reproduce the $\max |E|$ history up to $t \approx 9$ and $t\approx 14$ for $\bp = (0.01, 1)$ and $(0.03, 0.8)$, respectively. On the other hand, we observe that the TW-ROMs are able to reproduce the $\max |E|$ history for the entire time interval $\Tinterval$.

\begin{figure}[!ht]
    \centering
    \includegraphics[width=0.9\textwidth]{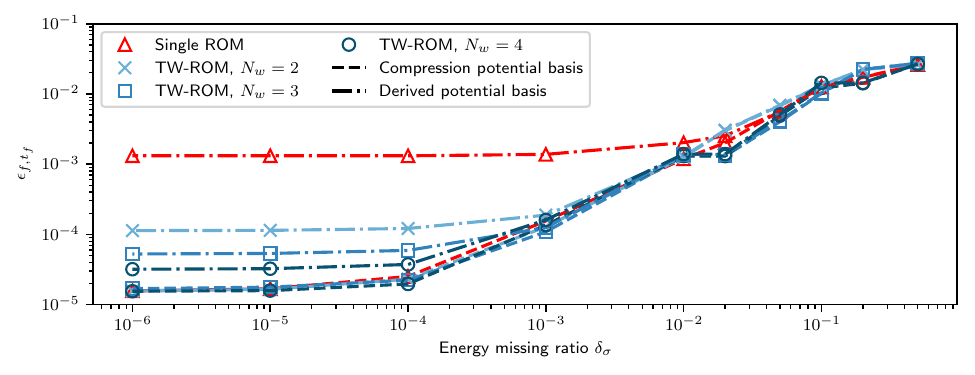}
    \caption{
    Reproduction problem of the Landau-damping problem with final time $t_f=15$ for parameter $\bp=(0.03,0.8)$. 
    Comparison of ROMs using two methods for constructing the potential reduced basis functions. The compression potential basis refers to the approach utilizing POD, while the derived potential basis involves solving Poisson problems.}
    \label{fig:ld_reprod_error_compare_basis}
\end{figure}

We further assess the performance of the single ROM and TW-ROMs in the reproduction problem for the parameter $\bp=(0.03,0.8)$ using the two
methods for constructing the potential reduced basis functions, as discussed in
Section~\ref{sec:solution_basis_construction}. In
Fig.~\ref{fig:ld_reprod_error_compare_basis}, the behavior of the $\epf$ with
respect to the $\delta_\sigma$ for the two approaches is shown. The compression
potential basis refers to the approach that utilizes POD, while the derived potential basis involves
solving Poisson problems. 
We find that for values of $\delta_\sigma  > \tento{-3}$, the performance of the ROMs is comparable for both.
For values of $\delta_\sigma \le \tento{-3}$, the ROMs with the compression potential basis are more accurate than the ROMs with the derived potential basis. In particular, we observe a significant difference in the performance of the single ROM. 

The differences in performance can be attributed to two factors. First, we find that the compression potential basis is more effective in representing the potential snapshots compared to the derived potential basis. In particular, for a given value of $n_\phi$, the projection error associated with the compression potential basis is five orders of magnitude smaller than that of the derived potential basis.
Second, we compare the singular values of the potential and distribution snapshot matrices. We find that for a given energy missing ratio, such as $\delta_\sigma = \tento{-4}$, (\ref{eq:energy_criteria}) indicates that the basis size for the potential is $n_\phi=3$, while the basis size for the distribution is $n_f=30$. 
With the derived potential basis, the reduced potential field uses
the same number of modes as the reduced distribution field. The inclusion of high-index modes can contribute to the error due to the nonlinear approximation in (\ref{eq:tensorial_approx}). This explains why the ROMs with derived potential basis have larger errors.

\subsubsection{Parametric problem}
We next consider the parametric problem, where the ROMs are deployed to predict solutions at parameters $\bp \in \mD$ other than the training parameters. 
To construct the ROMs, we consider $n_p = 4$ training parameters, namely
$\bp_1= (0.01,0.8)$, $\bp_2=(0.01,1.0)$, $\bp_3=(0.03,0.8)$, and
$\bp_4=(0.03,1.0)$. For each training parameter, $3000$ snapshots are
collected in the time interval $\Tinterval = (0, 15]$, leading to a total of $K=12000$ snapshots. The single ROM (\ref{eq:tensorial-rom}) is constructed by performing POD on the entire snapshot set.
For the TW-ROM (\ref{eq:TW--ROM}), 
the range of the physical time indicator $[\indicator_{\text{min}} = 0,~
\indicator_{\text{max}}=15)$ is partitioned uniformly into $N_w$ subintervals. 
$K=12000$ snapshots are then classified into $N_w$ groups based on the partition of the physical time indicator range, where the $j$-th group is 
\begin{align}
   \mG_j = \{(n,k) \in \mathbb{Z} \times \mathbb{N}(n_p): 0\le n \le N_t ~\text{and}~ t_n \in [\indicator_{j-1}, ~\indicator_j) \}.
\end{align}
Then, for $1 \le j \le N_w$, the TW-ROM is constructed by performing POD on the subset of snapshots whose indices belong to the group $\mG_j$.
In the online phase, 
the ROMs are deployed at the various parameters $\bp \in
\mD^h \subset \mD = [0.01,0.03] \times [0.8,1]$ with final time $t_f = 15$,
where $\mD^h$ is a discrete parameter space $\mD^h \subset \mD %=
%[0.01,0.03]\times[0.8,1]$ 
$ with $5$ and $5$ evenly distributed discrete points in
the respective parameter range, resulting in a total of $25$ parameters. For each parameter $\bp$, the corresponding reduced order model is used, where the end-point $T_j(\bp)$ of the temporal subinterval $\mT_j(\bp)$ is determined similarly as in the reproduction problem.

\begin{figure}[!ht]
    \centering
    \begin{subfigure}{0.49\columnwidth}
    \caption{Single ROM}
    \includegraphics[width=1\columnwidth]{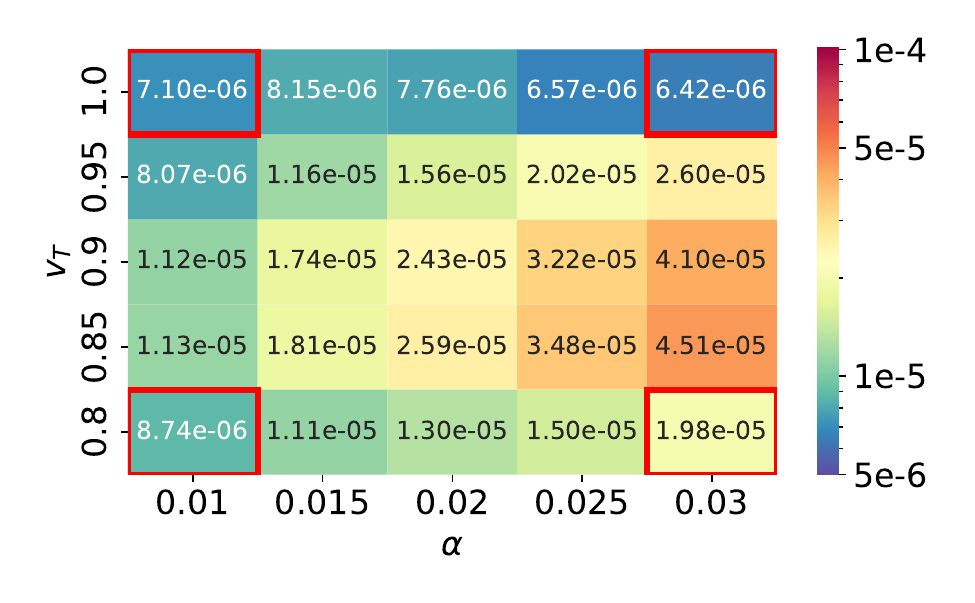}
    \end{subfigure}
    \begin{subfigure}{0.49\columnwidth}
    \caption{TW-ROM with $N_w=2$}
    \includegraphics[width=1\columnwidth]{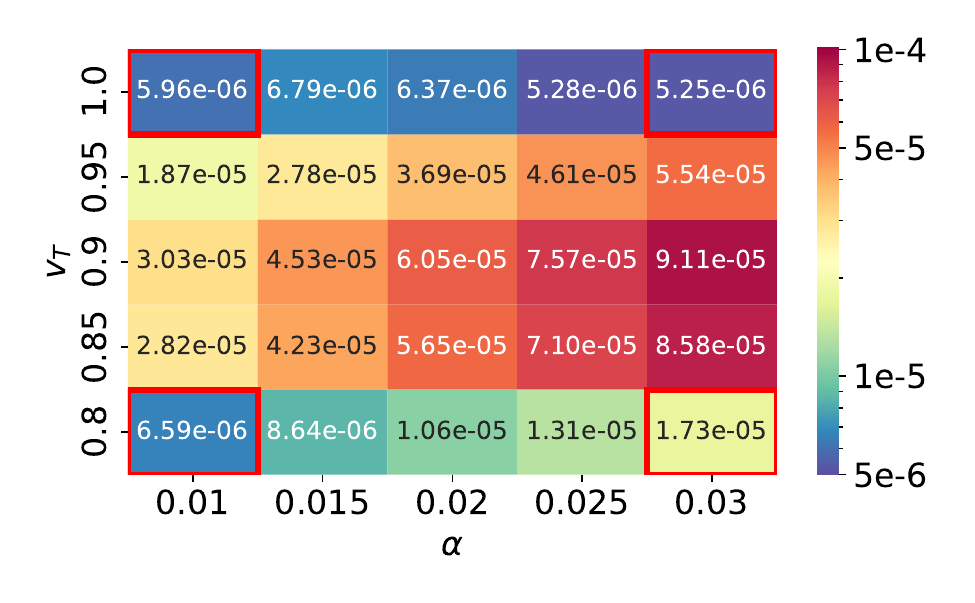}
    \end{subfigure}
    \begin{subfigure}{0.49\columnwidth}
    \caption{TW-ROM with $N_w=3$}
    \includegraphics[width=1\columnwidth]{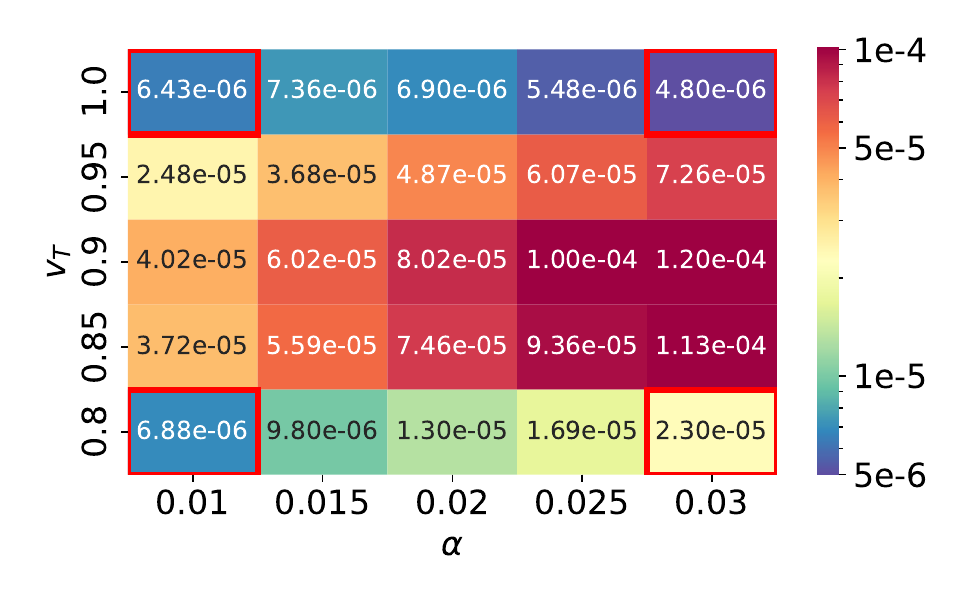}
    \end{subfigure}
    \begin{subfigure}{0.49\columnwidth}
    \caption{TW-ROM with $N_w=4$}
    \includegraphics[width=1\columnwidth]{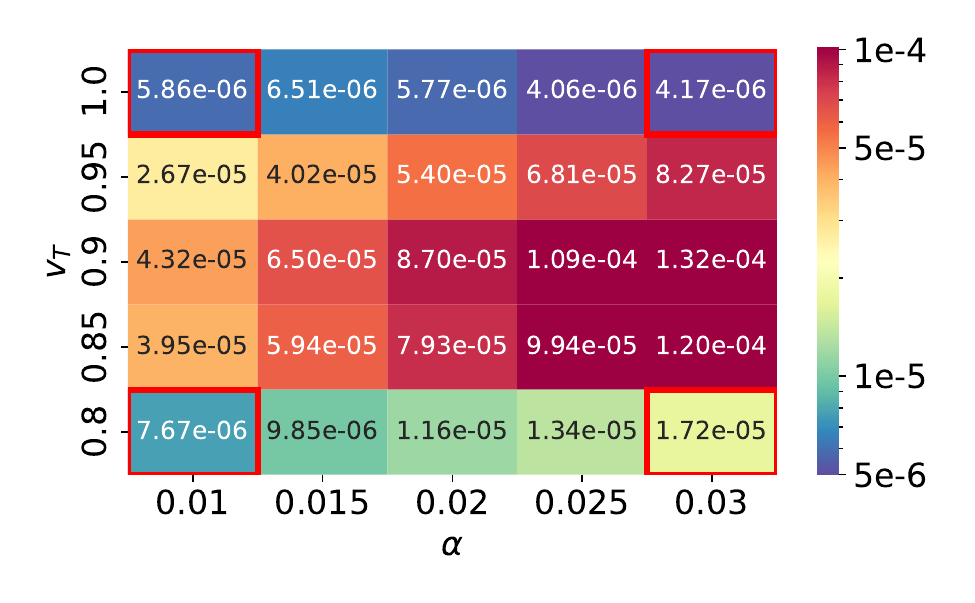}
    \end{subfigure} 
    \caption{Parametric problem of the Landau damping. The relative error in the distribution field at the final time $\epf$ of
    $\bp \in \mD^h$ for the single ROM and the TW-ROMs with various values of $N_w$.
    The energy missing ratio $\delta_\sigma$ is chosen to be $\tento{-4}$ for all the ROMs. The values of $\epf$ for the
    training parameters are enclosed with a red box.} 
    \label{fig:ld_pmor_error}
\end{figure}

Fig.~\ref{fig:ld_pmor_error} displays the relative error in the distribution
field at the final time $\epf$ of $25$ parameters $\bp \in \mathcal{D}^h$ for
the single ROM and the TW-ROMs with three values of $N_w$. The energy missing ratio
$\delta_\sigma$ is chosen to be $\tento{-4}$ for all the ROMs. 
We observe that the TW-ROMs have smaller errors at the four training parameters 
and testing parameters, $\bp \in \{0.015,0.02,0.025\} \times \{0.8,1.0\}\}$, compared to the single ROM. 
In addition, by refining a fixed temporal partition—--specifically, by using values of $N_w=2$ and $N_w=4$—--we find that the TW-ROM with $N_w=4$ performs slightly better than the TW-ROM with $N_w=2$ at the training parameters.
For the testing parameters, $\bp \in \{0.01,0.015,0.02,0.025, 0.03\} \times \{0.85, 0.9, 0.95\}$, the single ROM has smaller errors. Additionally, increasing the value of $N_w$ slightly deteriorates the performance of the TW-ROM.
Overall, we observe that both the single ROM and the TW-ROMs achieve a maximum error of less than $0.01\%$ for parameters $\bp \in \mD^h$.

We further calculate the average speed-up for the 25 parameters in $\mathcal{D}^h$ for each ROM. We find that the single ROM achieves a speed-up of 20, while the TW-ROMs demonstrate improved speed-ups of 36, 46, and 54 times for \(N_w = 2\), \(N_w = 3\), and \(N_w = 4\), respectively.

\begin{figure}[!ht]
    \centering
    \begin{subfigure}{0.49\columnwidth}
    \caption{$\bp = (0.03,0.8)$}
    \includegraphics[width=1\columnwidth]{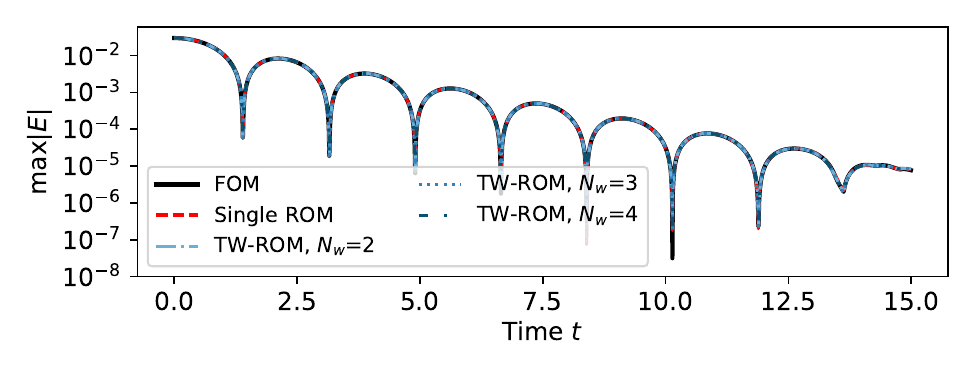}
    \end{subfigure}
    \begin{subfigure}{0.49\columnwidth}
    \caption{$\bp = (0.01, 1)$}
    \includegraphics[width=1\columnwidth]{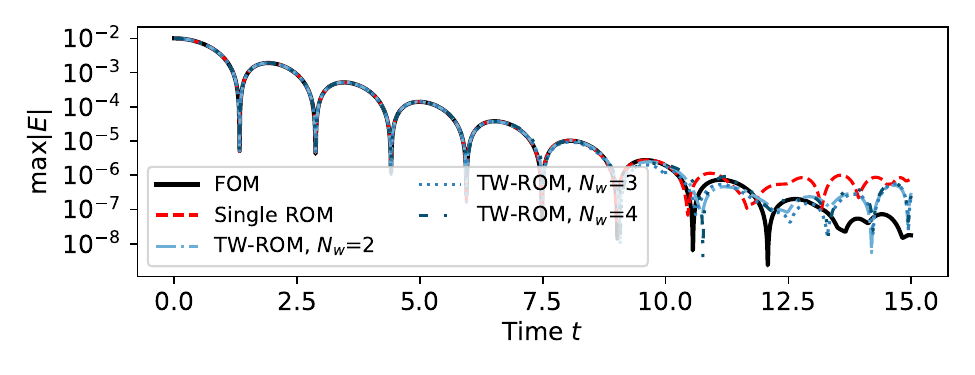}
    \end{subfigure}
    \begin{subfigure}{0.49\columnwidth}
    \caption{$\bp = (0.03, 0.9)$}
    \includegraphics[width=1\columnwidth]{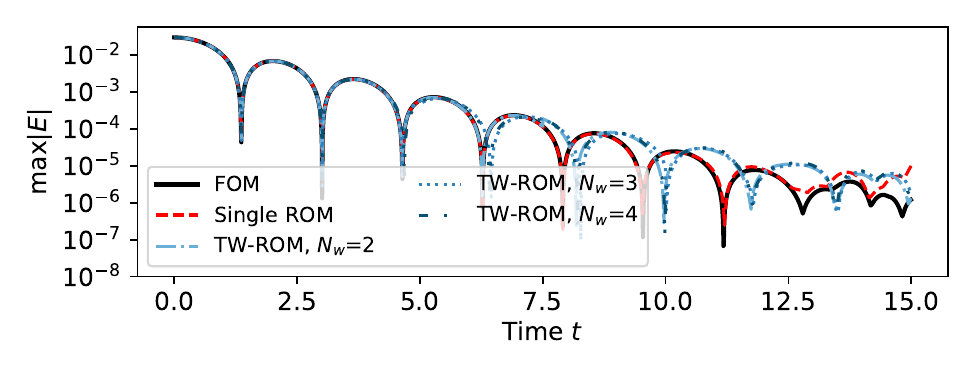}
    \end{subfigure}
    \begin{subfigure}{0.49\columnwidth}
    \caption{$\bp = (0.03, 0.85)$}
    \includegraphics[width=1\columnwidth]{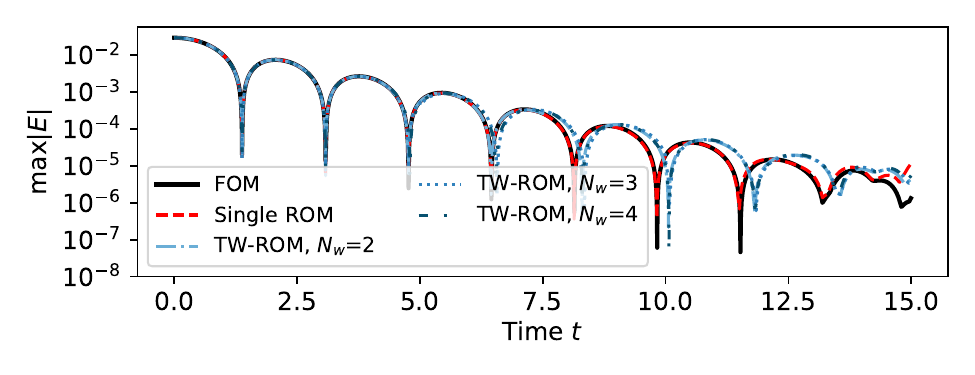}
    \end{subfigure}
    \caption{Parametric problem of the Landau damping. $\max |E|$ behavior in the time interval $\Tinterval$ of the single ROM
    (\ref{eq:1D--1V-vlasov-rom}) and the TW-ROM (\ref{eq:TW--ROM}) with the energy missing ratio $\delta_\sigma =\tento{-4}$ for the two training parameters $\bp=(0.03,0.8)$ and $(0.01,1)$, and the two testing parameters $\bp=(0.03,0.9)$ and $(0.03,0.85)$.
    } %\pht{TODO: Update plots with upw5 results.}}
    \label{fig:ld_pmor_maxe}
\end{figure}

With $\delta_\sigma = \tento{-4}$, we further examine the approximated $\max |E|$
in the time interval $\Tinterval$ for the two training parameters
$\bptrain=(0.03,0.8)$ and $\bptrain = (0.01,1)$ and the two testing parameters
$\bptest=(0.03,0.9)$ and $\bptest=(0.03,0.85)$ in Fig.~\ref{fig:ld_pmor_maxe}.
The two testing parameters are specifically chosen because they correspond to ROMs with large $\epf$ values. 
We find that both the single ROM and the TW-ROMs can reproduce the $\max |E|$ history for $\bptrain=(0.03,0.8)$. For $\bptrain=(0.01,1)$, $\max |E|$ of all the ROMs deviate away from the FOM $\max |E|$ after $t \approx 10$ and a much larger discrepancy is observed in the single ROM.
For the two testing parameters $\bptest$, both the single ROM and the TW-ROMs struggle to accurately predict the behavior of $\max |E|$. We observe that the $\max |E|$ of the single ROM grows towards the end of $\Tinterval$. Although $\max |E|$ of the TW-ROMs remains stable, there is a shift in the period of $\max |E|$.
While there are discrepancies in the $\max |E|$ behavior compared to the FOM, the decay rate of $\max |E|$ predicted by the ROMs is similar to that of the FOM. 

\subsection{2D (1D1V) Vlasov\dash Poisson equations - Two-stream instability}
\label{ssec:two_stream}

The formation of the two-stream instability occurs due to the nonlinear
evolution of (\ref{eq:1D--1V-vlasov}). This instability is a well-known
phenomenon in plasma physics, which is generated by two counter-streaming beams.
In this process, the kinetic energy of particles excites a plasma wave, which
then transfers to electrostatic potential energy~\cite{anderson2001tutorial}. In this study, we focus on the temporal 
interval that includes three stages, namely, short
transient, growth, and statistically stationary stages, which require many
reduced bases for ROMs to accurately capture the behavior in each stage~\cite{hesthaven2023adaptive}. The initial condition is 
\begin{align}
   f(x,v) = \frac{n_0}{\sqrt{2\pi \vT^2}}(1 + \alpha \cos(2 k \pi \frac{x}{L}) )\left[ \exp(-\frac{(v-v_0)^2}{2 \vT^2}) + \exp(-\frac{(v+v_0)^2}{2 \vT^2}) \right], \label{eq:sc_initial}
\end{align}
where $L=2\pi$, $n_0=8$ and $v_0=1$ is the initial velocity displacement in phase space. The wavenumber $k$ of the perturbation is set to $1$. The perturbation amplitude $\alpha$ and the thermal velocity $\vT$ are the study parameters $\bp=(\alpha, \vT)$, with $\bp \in \mD = [0.001,0.0025] \times [0.08,0.1]$.

We consider the computational domain $(x,v) \in \Omega \coloneqq [0,2\pi]\times [-3.5\vT, 3.5\vT]$ with periodic boundary conditions in the $x$-direction and homogeneous Dirichlet boundary conditions in the $v$-direction. The FOM has $N_f = 256 \times 256 = 65536$ degrees of freedoms. A uniform time step $\Delta t = 0.0025$ is used for the evolution of the particle's distribution over the time interval $\Tinterval = (0, t_f]$ with $t_f= 10$. 

\begin{figure}[!ht]
    \centering
    \includegraphics[width=0.9\textwidth]{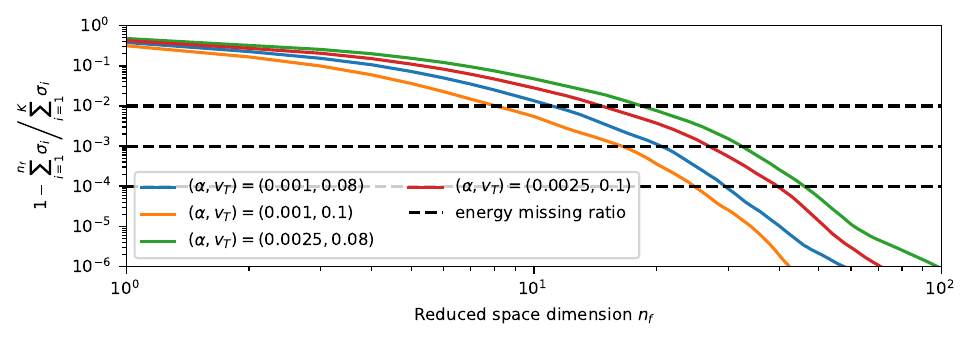}
    \caption{Two-stream instability. The behavior of 
    $1-\sum^{n_f}_{i=1} \sigma_i /\sum^K_{i=1} \sigma_i$ as a function of the reduced space dimension $n_f$ for four parameters $\bp=(0.001,0.08),(0.001,0.1),(0.0025,0.08)$ and $(0.0025,0.1)$.}
    \label{fig:selfconsistent_singular_reprod1}
\end{figure}

In Fig.~\ref{fig:selfconsistent_singular_reprod1}, we investigate the reducibility of
the problem for four parameters, namely, $\bp = (0.001,0.08)$, $(0.001,0.1)$, $(0.0025,0.08)$ and $(0.0025,0.1)$. For each
parameter $\bp$, we compute the singular values of the snapshots matrix $\bF \in
\mathbb{R}^{N_f \times K}$, where $K=4000$ snapshots are collected within
$\Tinterval$. The behavior of $1-\sum^{n_f}_{i=1}\sigma_i/\sum^K_{i=1}\sigma_i$
with respect to the reduced space dimension $n_f$ is shown with three values of the energy
missing ratios, that is, $\delta_\sigma = \tento{-2},~\tento{-3}$ and
$\tento{-4}$, represented by the red dashed line. 
The results indicate that the global representation is less efficient when the perturbation amplitude $\alpha$ is larger and the thermal velocity $v_T$ is smaller. Specifically, this is in comparison to the case where $\bp = (0.001, 0.1)$.
This is expected because, within a fixed time interval $\Tinterval$, a larger value of $\alpha$ tends to result in the snapshot matrix containing a greater number of highly nonlinear solutions from the statistically stationary stage. Additionally, when $v_T$ is small, the solutions that appear in the matrix are generally less smooth.

We note that from the evolution of the electric field energy $\int E^2~dx$ shown
in Fig.~\ref{fig:fom_maxE}, the time $t$ required to reach a statistically
stationary stage varies by a factor of $1$ across different parameter
variations. Additionally, the reducibility illustrated in
Fig.~\ref{fig:selfconsistent_singular_reprod1} indicates a degree of variability
in the solutions within the considered parameter range. 

\subsubsection{Reproduction problem}
\begin{figure}[!ht]
    \centering
    \begin{subfigure}{0.20\columnwidth}
    \caption{$t=3$}
    \includegraphics[width=1\textwidth]{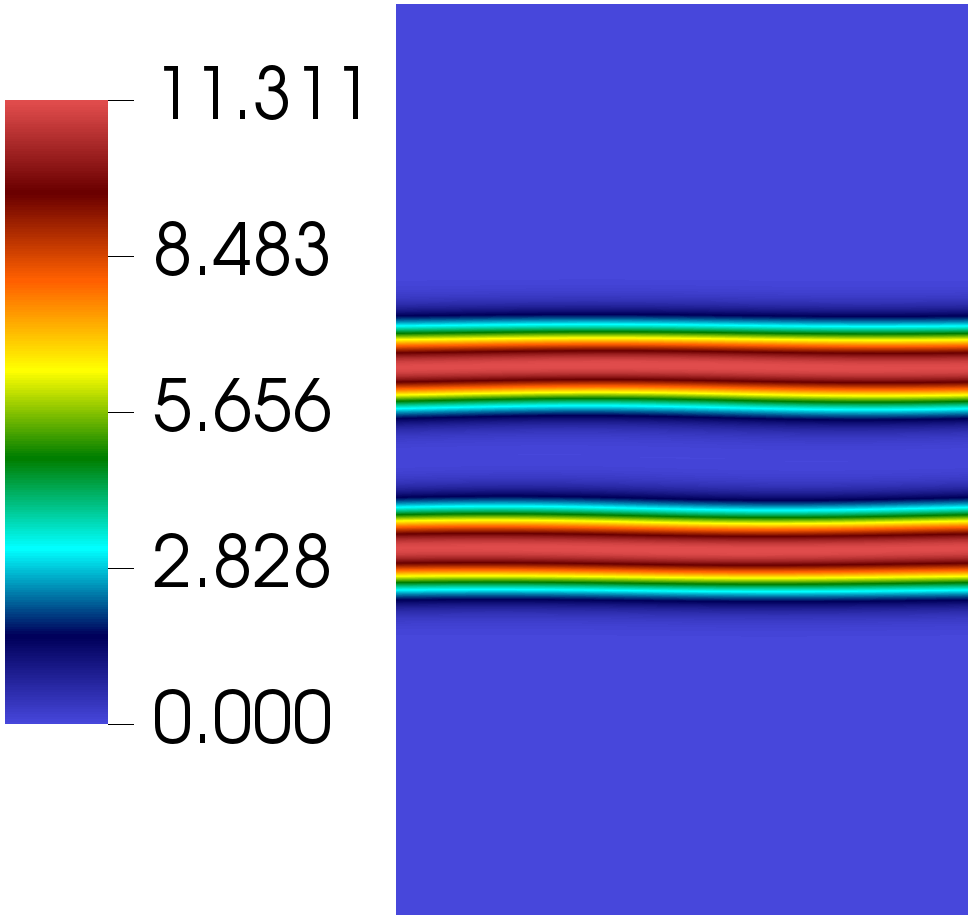}
    \end{subfigure}
    \begin{subfigure}{0.20\columnwidth}
    \caption{$t=4$}
    \includegraphics[width=1\textwidth]{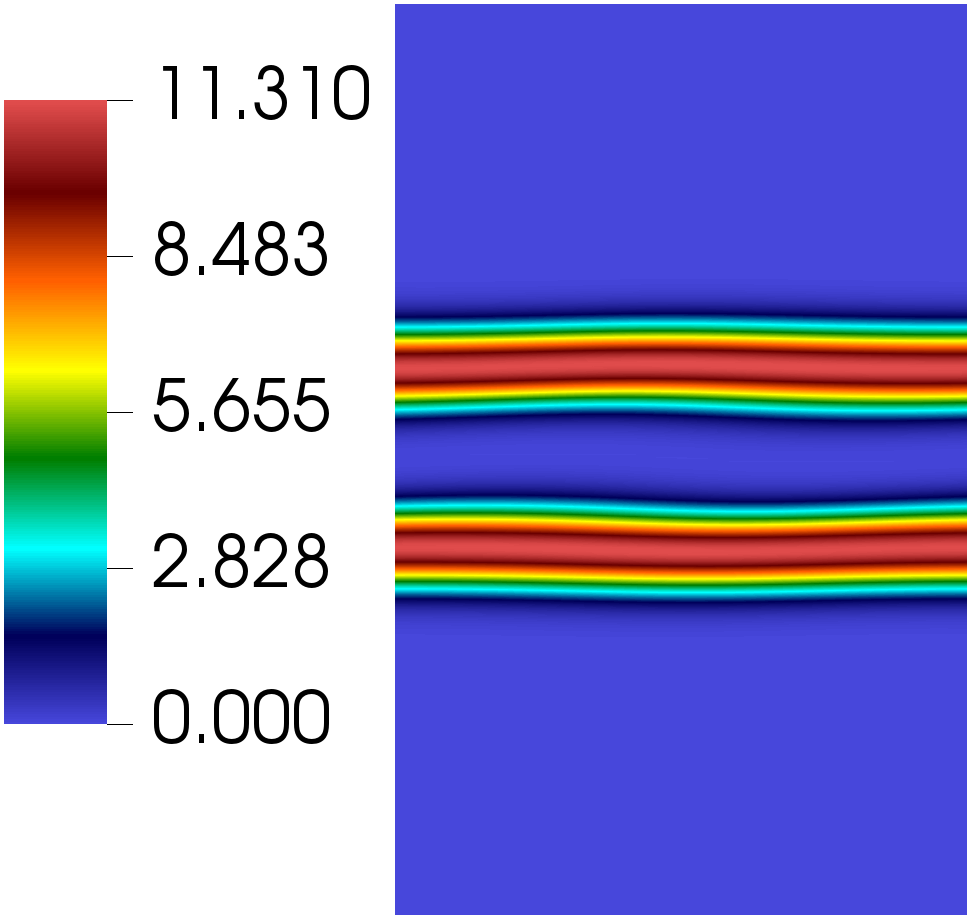}
    \end{subfigure}
    \begin{subfigure}{0.20\columnwidth}
    \caption{$t=5$}
    \includegraphics[width=1\textwidth]{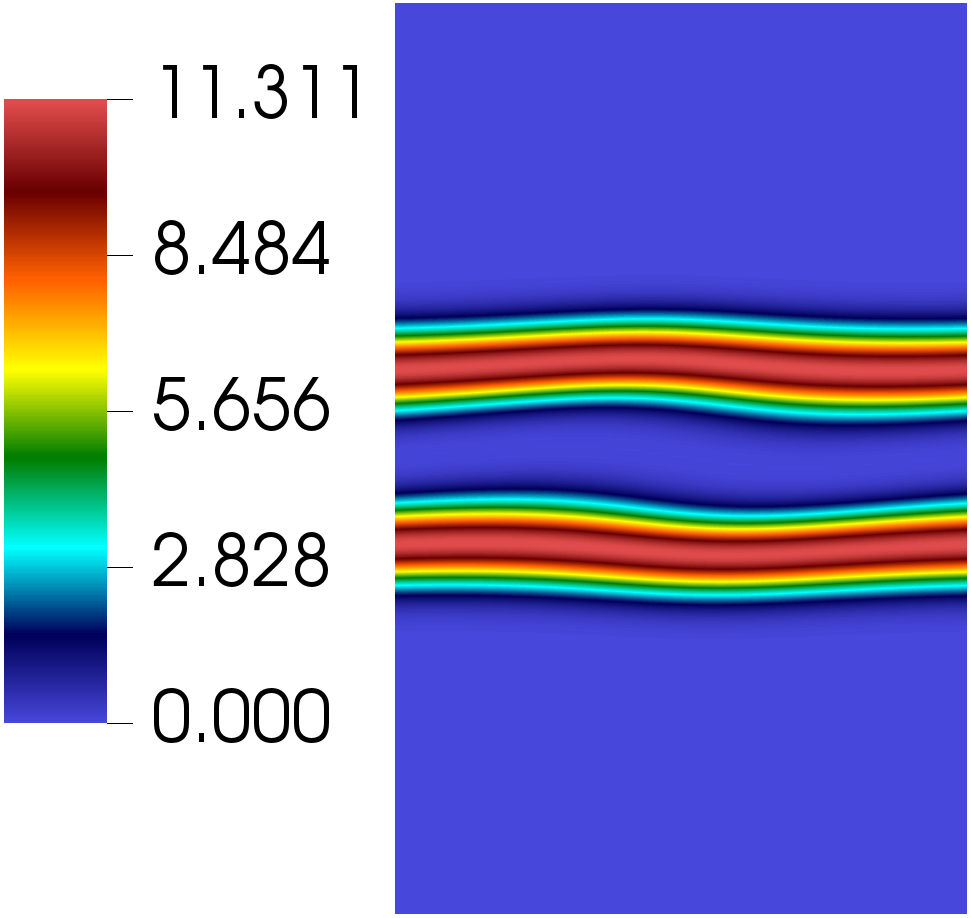}
    \end{subfigure}
    \begin{subfigure}{0.20\columnwidth}
    \caption{$t=6$}
    \includegraphics[width=1\textwidth]{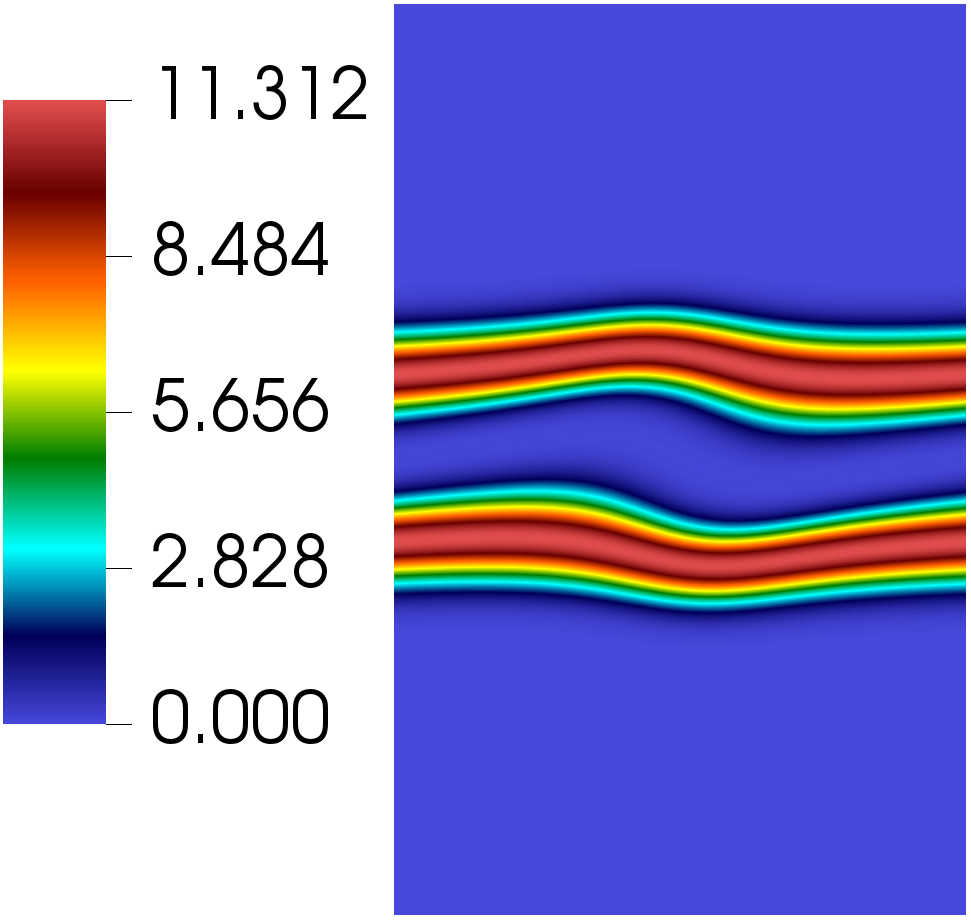}
    \end{subfigure}\\
    \begin{subfigure}{0.20\columnwidth}
    \caption{$t=7$}
    \includegraphics[width=1\textwidth]{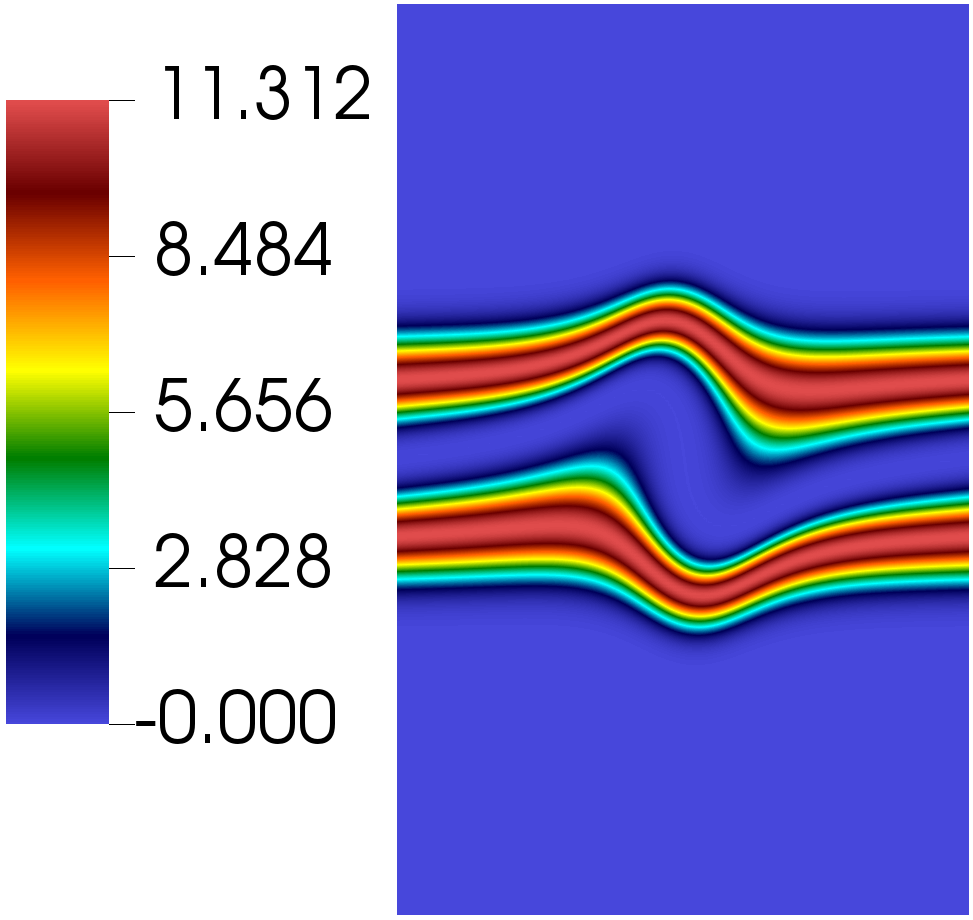}
    \end{subfigure}
    \begin{subfigure}{0.20\columnwidth}
    \caption{$t=8$}
    \includegraphics[width=1\textwidth]{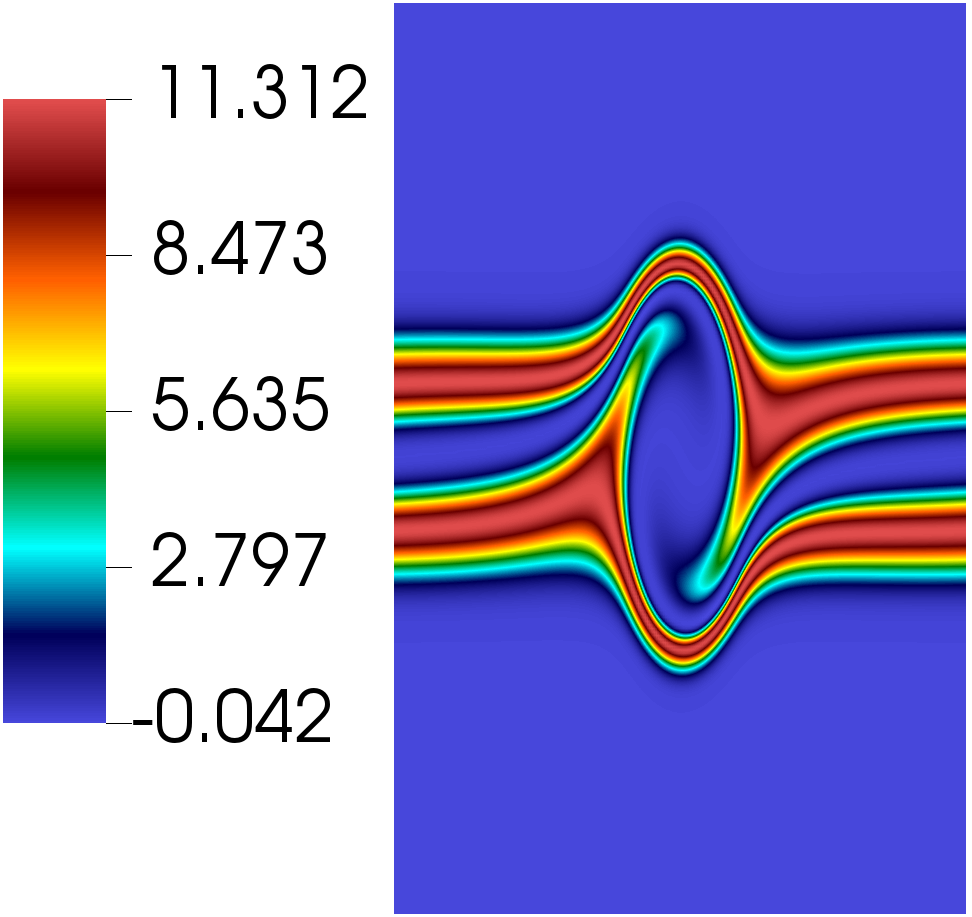}
    \end{subfigure}
    \begin{subfigure}{0.20\columnwidth}
    \caption{$t=9$}
    \includegraphics[width=1\textwidth]{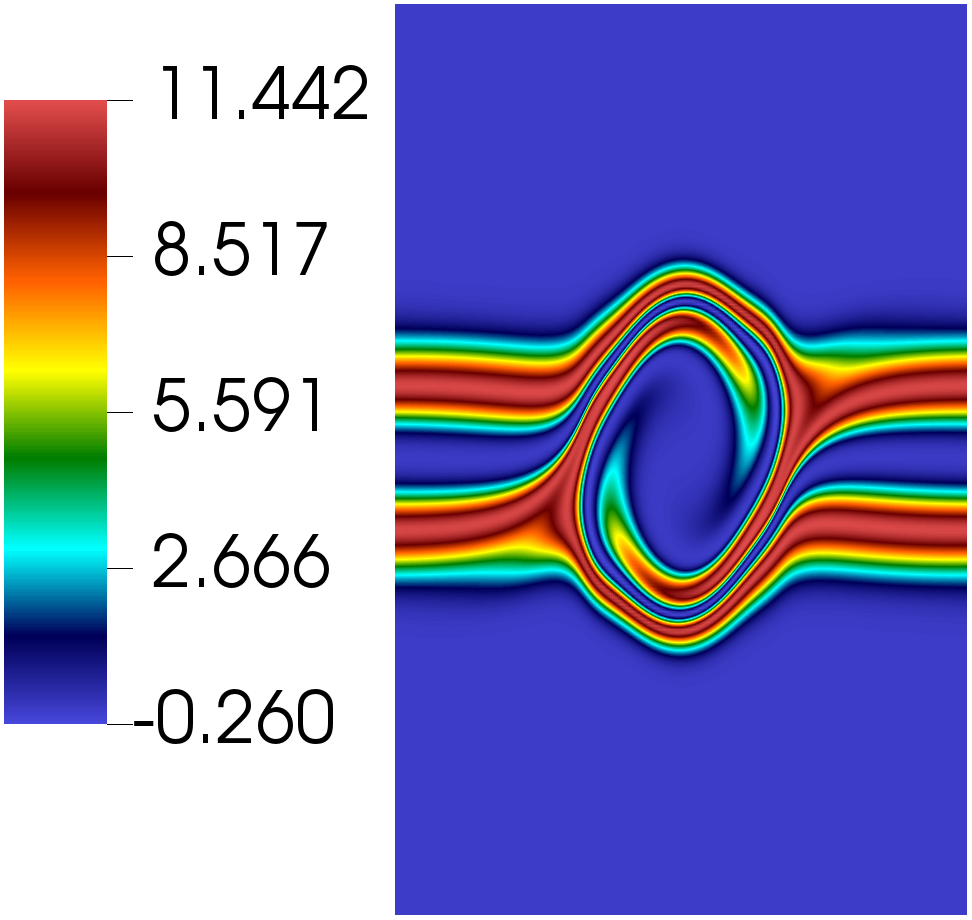}
    \end{subfigure}
    \begin{subfigure}{0.20\columnwidth}
    \caption{$t=10$}
    \includegraphics[width=1\textwidth]{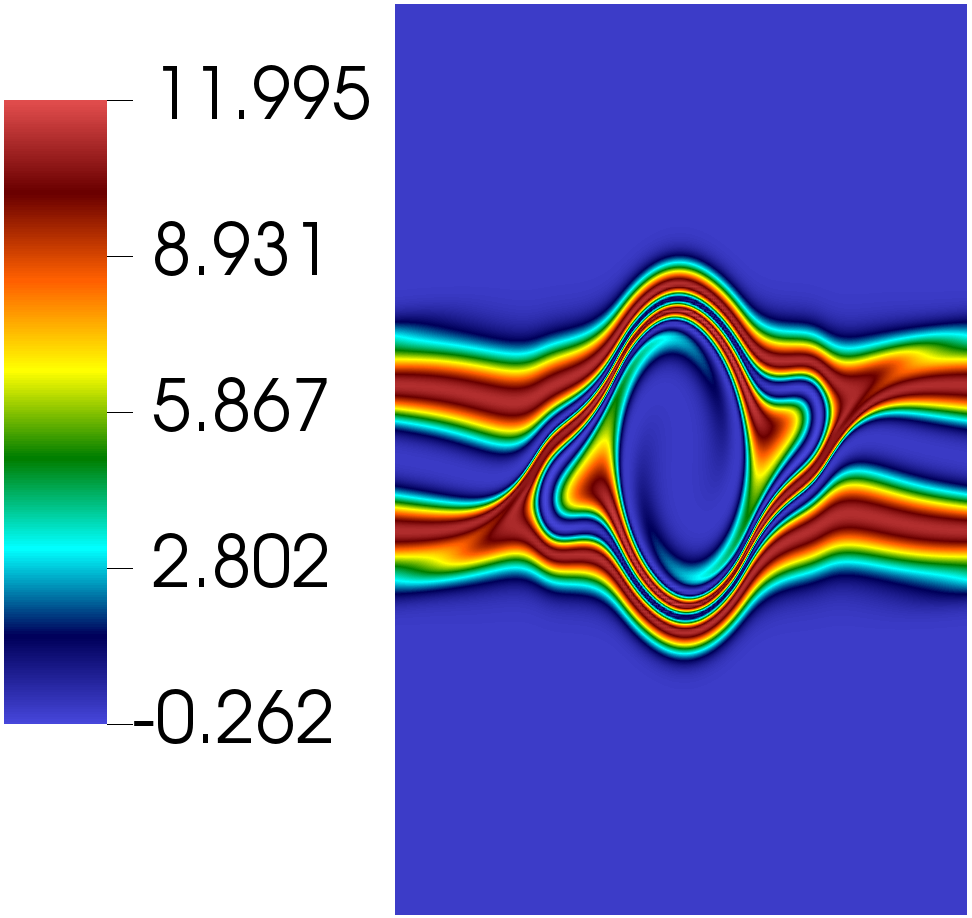}
    \end{subfigure}
    \caption{Two-stream instability with $\alpha=0.0025$ and $\vT=0.08$. FOM distribution field at different times $t$.}
    \label{fig:sc_reprod1_snapshots}
\end{figure}

We consider the reproduction problem at $\bp=(0.0025, 0.8)$. 
The FOM distribution field at time $t=3,4,\ldots,9,10$ is shown in  Fig.~\ref{fig:sc_reprod1_snapshots}. We observe solutions in three stages, namely, the short transient stage ($t=3,4,5$), the growth stage ($t=6,\ldots,8$), and the statistically stationary stage ($t=9,10$). 

The single ROM is constructed by performing POD on $K=4000$ snapshots collected in the time interval $\Tinterval$. 
For the TW-ROM (\ref{eq:TW--ROM}), the range of the 
physical time indicator $[\indicator_{\text{min}} = 0,~
\indicator_{\text{max}}=10)$ is partitioned uniformly into $N_w$ subintervals. 
$K=4000$ snapshots are then classified into groups based on the partition of the indicator range. 
Then, for $1 \le j \le N_w$, the TW-ROM is constructed by performing POD on the subset of snapshots whose indices belong to the group $\mG_j$. In the online phase, the corresponding reduced order model is used, where the end-point $T_j(\bp)$ of the temporal subinterval $\mT_j(\bp)$ is defined as the time instance $\widetilde{t_n}(\bp)$ when $\indicator(\tbff_n(\bp_k), \widetilde{t_n}, \bp_k)$ first exceeds $\indicator_j$, at which we increment to the next subinterval $\mT_{j+1}(\bp)$.

For the EW-ROM (\ref{eq:EW--ROM}), $K=4000$ snapshots are classified into $N_w=5$ groups defined in Section~\ref{ssec:EW--ROM}. 
Then, for $1 \le j \le 5$, the EW-ROM is constructed by performing POD on the subset of snapshots whose indices belong to the group $\mG_j$. 
In the online phase, the corresponding reduced order model is used, where 
the end-point $T_2(\bp)$ of the temporal subinterval $\mT_2(\bp)$ is defined as the time instance $\widetilde{t_n}(\bp)$ when the electric field energy $\indicator(\tbff_n(\bp_k), \widetilde{t_n}, \bp_k)$ first exceeds $1$, at which we increment to the next subinterval $\mT_{3}(\bp)$. 
The end-point $T_3(\bp)$ of the temporal subinterval $\mT_3(\bp)$ is defined as the time instance $\widetilde{t_n}(\bp)$ when the slope of the electric field energy $\partial \indicator(\tbff_n(\bp_k), \widetilde{t_n}, \bp_k)/\partial t$ becomes less than $0$, at which we increment to the next subinterval $\mT_{4}(\bp)$.
The end-point $T_4(\bp)$ of the temporal subinterval $\mT_4(\bp)$ is defined as the time instance $\widetilde{t_n}(\bp)$ when the slope of the electric field energy $\partial \indicator(\tbff_n(\bp_k), \widetilde{t_n}, \bp_k)/\partial t$ becomes greater than $0$, at which we increment to the final subinterval $\mT_{5}(\bp)$.
\begin{figure}[!ht]
    \centering
    \begin{subfigure}{0.9\columnwidth}
    \centering
    \caption{}
    \includegraphics[width=1\textwidth]{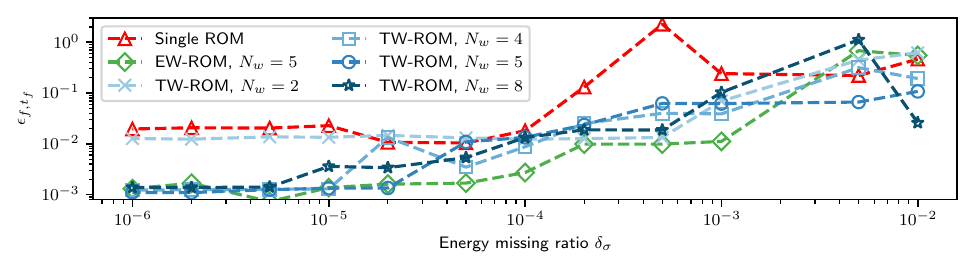}
    \end{subfigure}
    \begin{subfigure}{0.9\columnwidth}
    \centering
    \caption{}
    \includegraphics[width=1\textwidth]{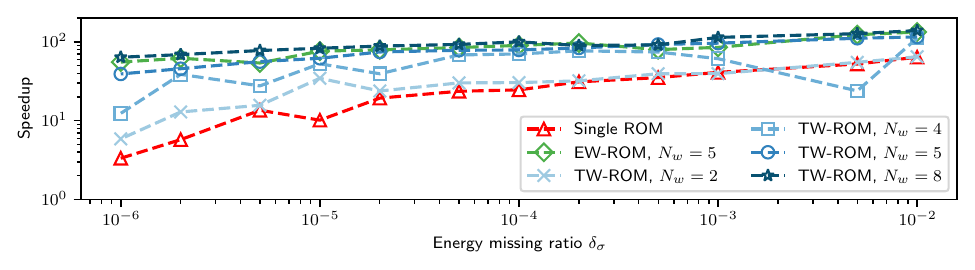}
    \end{subfigure}
    \caption{Reproduction problem of the two-stream instability for parameter $\bp = (\alpha,v_T) = (0.0025,0.08)$. (a) The behavior of the relative error in the distribution field at the final time $\epf$ and (b) the speed-up as a function of the energy missing ratio $\delta_\sigma$ for the single ROM, the TW-ROMs and the EW-ROM.} 
    \label{fig:sc_error_compare}
\end{figure}

In Fig.~\ref{fig:sc_error_compare} (a), the behavior of the relative error $\epf$
with respect to $\delta_\sigma$ for the ROMs is shown. 
The results show that the $\epf$ for the single ROM oscillates as $\delta_\sigma$ decreases and can achieve an error of $2\%$ with $\delta_\sigma = 10^{-5}$. In comparison, the results for the TW-ROMs show improved performance over the single ROM for the range $\tento{-3}\le \delta_\sigma \le \tento{-4}$. For $\delta_\sigma \le \tento{-4}$, the TW-ROM with $N_w=2$ achieves an accuracy similar to the single ROM whereas the TW-ROMs with $N_w>2$ are able to achieve a better accuracy. 
For the EW-ROM, we find it outperforms the TW-ROM with the same or larger $N_w$ values for $\tento{-3} \le \delta_\sigma \le \tento{-5}$. This demonstrates the effectiveness of using the electric field energy indicator to classify the FOM states.

In Fig.~\ref{fig:sc_error_compare}
(b), the behavior of the speed-up with respect to $\delta_\sigma$ for the ROMs is shown. The results show that the single ROM has the least speed-up. In contrast, the speed-up of the TW-ROM improves as the number of subintervals $N_w$ increases. Additionally, the results show that the speed-up of the EW-ROM is similar to the TW-ROM with the same ($N_w=5$) or larger ($N_w=8$) number of subintervals for $\delta_\sigma \ge \tento{-5}$. However, for $\delta_\sigma \le \tento{-5}$, the EW-ROM has a better speed-up compared to the TW-ROM with $N_w=5$.

\begin{figure}[!ht]
    \begin{subfigure}{0.68\columnwidth}
    \centering
    \caption{} 
    \includegraphics[width=1\textwidth]{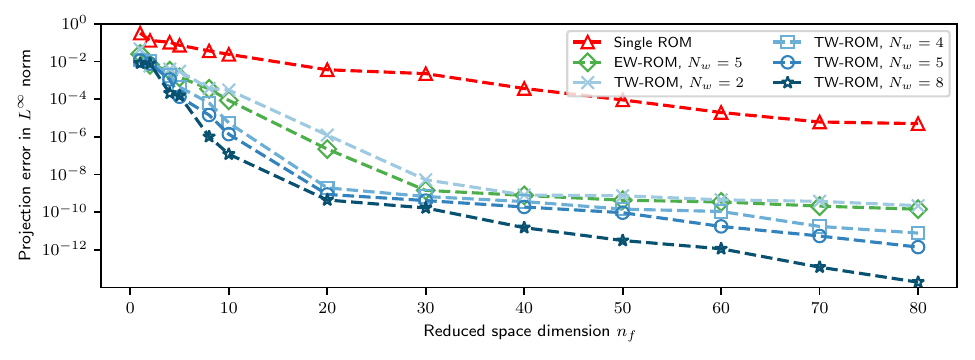}
    \end{subfigure}
    \begin{subfigure}{0.31\columnwidth}
    \centering
    \caption{}
    \vspace{0.26cm}
    \includegraphics[width=0.8\textwidth]{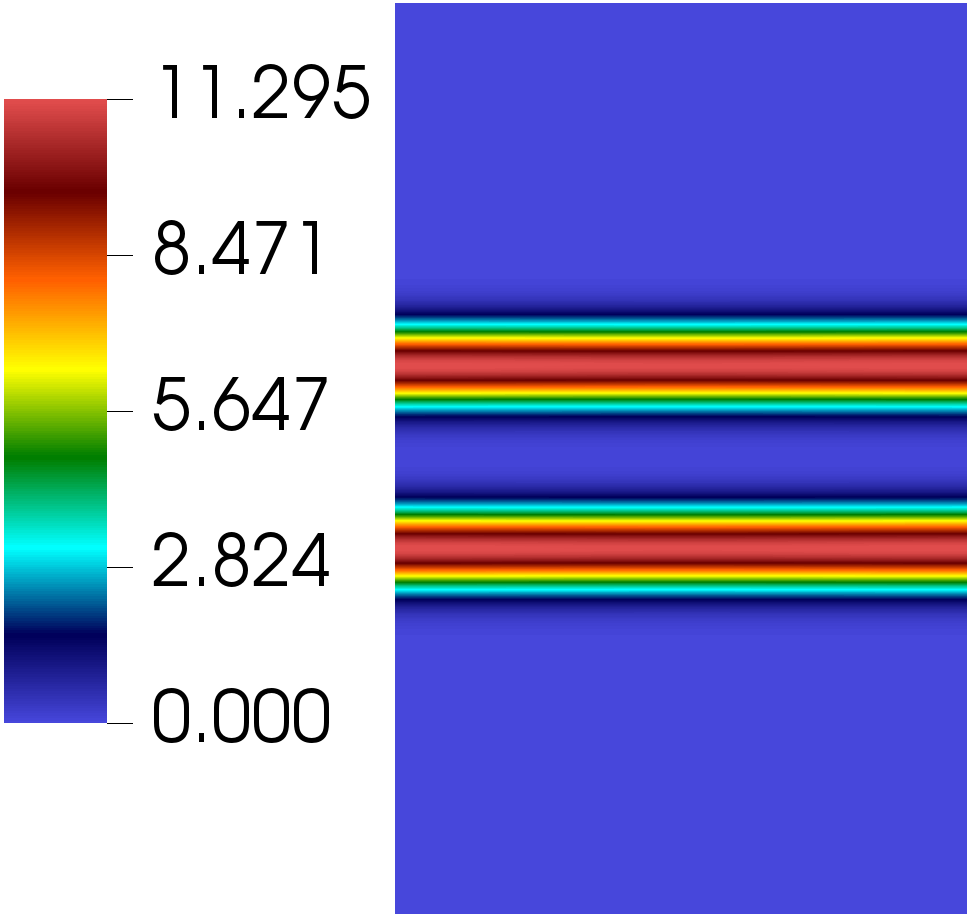}
    \end{subfigure}
    \caption{Two-stream instability with $\bp = (\alpha,v_T) = (0.0025,0.08)$. (a) The behavior of the projection error of the initial distribution state in $L^\infty$ norm with respect to reduced space dimension $n_f$. (b) Initial distribution state. } 
    \label{fig:sc_init_coef_reprod1}
\end{figure}

We further investigate the approximation of the reduced space of the single ROM, the TW-ROMs, and the EW-ROM to the initial distribution state (\ref{eq:sc_initial}). The relative error between the true initial state and the projected initial state onto the reduced space, denoted as the projection error, is shown as a function of the reduced space dimension $n_f$ for all the ROMs in Fig.~\ref{fig:sc_init_coef_reprod1} (a).
We observe that the projection error decreases as the reduced space dimension increases in all the ROMs. However, the projection error with the global reduced space (single ROM) saturates at a level that is four orders of magnitude larger compared to the EW-ROM and the TW-ROMs. This shows the global reduced space is less effective in constructing the simple initial distribution state shown in Fig.~\ref{fig:sc_init_coef_reprod1} (b). For the TW-ROMs and the EW-ROM, because of the temporal domain decomposition, the reduced space that is used to approximate the initial state is from the 1st group $\mG_1$, which does not contain solutions in the growth and statistically stationary stages. This makes the representation more effective and leads to a smaller projection error. Additionally, we find improvement in the projection error as the number of subintervals $N_w$ increases. This is because the snapshots in the 1st group are more correlated to the initial state as the number of subintervals $N_w$ increases. We note that the EW-ROM has a projection error similar to the TW-ROM with $N_w=2$ but not with $N_w=5$. This difference arises from the classification described in Section~\ref{ssec:EW--ROM}. Specifically, the first group $\mG_1$ contains snapshots from the time interval $[0, 4]$ which is similar to the first group of the TW-ROM with $N_w=2$.

\begin{figure}[!ht]
    \centering
    \includegraphics[width=0.9\textwidth]{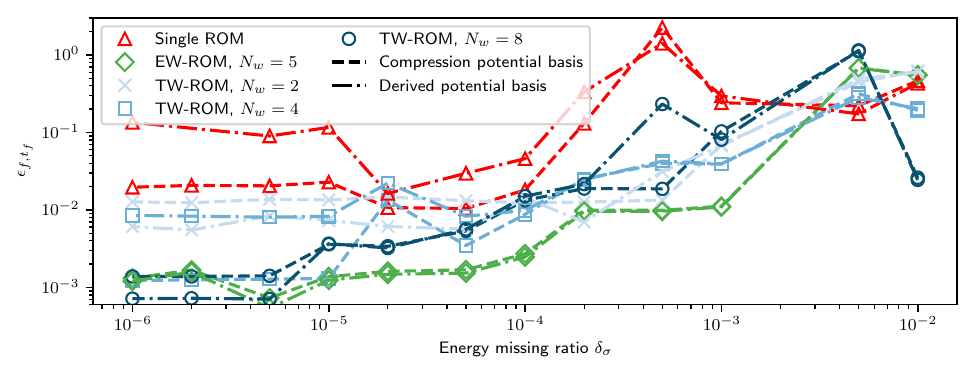}
    \caption{
    Reproduction problems of the two-stream instability with final time $t_f=15$ and parameter $\bp=(0.03,0.8)$. 
    Comparison of the ROMs using two methods for constructing potential reduced basis functions. The compression potential basis refers to the approach utilizing POD, while the derived potential basis involves solving Poisson problems.}
    \label{fig:sc_reprod_error_compare_basis}
\end{figure}

Similar to the Landau-damping problem, we further assess the
performance of ROMs using the two methods for constructing the potential reduced basis
functions, as discussed in Section~\ref{sec:solution_basis_construction}.
In Fig.~\ref{fig:sc_reprod_error_compare_basis}, the behavior of the $\epf$ with respect to $\delta_\sigma$ for the two methods is shown. We find that for values of $\delta_\sigma  > 5\times \tento{-3}$, 
the performance of the ROMs is comparable for both. For values of $\delta_\sigma \le 5 \times \tento{-3}$, 
%the small number of windows, that is, single ROM and TW-ROM with $N_w=2,4$, 
we find that the single ROM and TW-ROMs with $N_w=2,4$, using the compression potential basis, are more accurate than the ROMs with the derived potential basis. In particular, we observe a significant difference in the performance of the single ROM. 
For the EW-ROM and the TW-ROM with $N_w=8$, there is no significant difference between the two approaches. 
As discussed in the Landau-damping problem, the difference in performance can be attributed to two factors. First, for a given $n_\phi$ value, the projection error associated with the compression potential basis is two orders of magnitude smaller than that of the derived potential basis. Second, from the singular values of the potential and distribution snapshot matrices, we find that for $\delta_\sigma = \tento{-4}$, (\ref{eq:energy_criteria}) indicates that the basis size for the potential is $n_\phi=16$, while the basis size for the distribution is $n_f=55$. 
The inclusion of high-index modes can contribute to the error due to the nonlinear approximation in (\ref{eq:tensorial_approx}). 

\subsubsection{Parametric problem}
The main interest of parametric studies in the two-stream instability is to investigate how the perturbation amplitude $\alpha$ and thermal velocity $\vT$ affect the solutions and the growth rate of the instability. 

To construct the ROMs, we consider $n_p = 4$ training parameters, namely
$\bp_1= (0.001,0.08)$, $\bp_2=(0.0025,0.08)$, $\bp_3=(0.001,0.1)$, and
$\bp_4=(0.0025,0.1)$ and for each training parameter, $4000$ snapshots are
collected in the time interval $\Tinterval = (0, 10]$, leading to a total of $K=16000$ snapshots. The single ROM (\ref{eq:tensorial-rom}) is constructed by performing POD on the entire snapshot set.
For the TW-ROM (\ref{eq:TW--ROM}), 
the range of the physical time indicator $[\indicator_{\text{min}} = 0,~
\indicator_{\text{max}}=10)$ is partitioned uniformly into $N_w$ subintervals. 
$K=16000$ snapshots are then classified into $N_w$ groups based on the partition of the physical time indicator range, where the $j$-th group is 
\begin{align}
   \mG_j = \{(n,k) \in \mathbb{Z} \times \mathbb{N}(n_p): 0\le n \le N_t ~\text{and}~ t_n \in [\indicator_{j-1}, ~\indicator_j) \}.
\end{align}
Then, for $1 \le j \le N_w$, the TW-ROM is constructed by performing POD on the subset of snapshots whose indices belong to the group $\mG_j$. For the EW-ROM (\ref{eq:EW--ROM}), 
$K=16000$ snapshots are classified into $N_w=5$ groups defined in Section~\ref{ssec:EW--ROM}.
Then, for $1 \le j \le N_w$, the TW-ROM and the EW-ROM are constructed by performing POD on the subset of snapshots whose indices belong to the group $\mG_j$. 
\begin{figure}[!ht]
    \centering
    \begin{subfigure}{0.49\columnwidth}
    \caption{Single ROM, $\delta_\sigma = \tento{-3}$}
    \includegraphics[width=1\columnwidth]{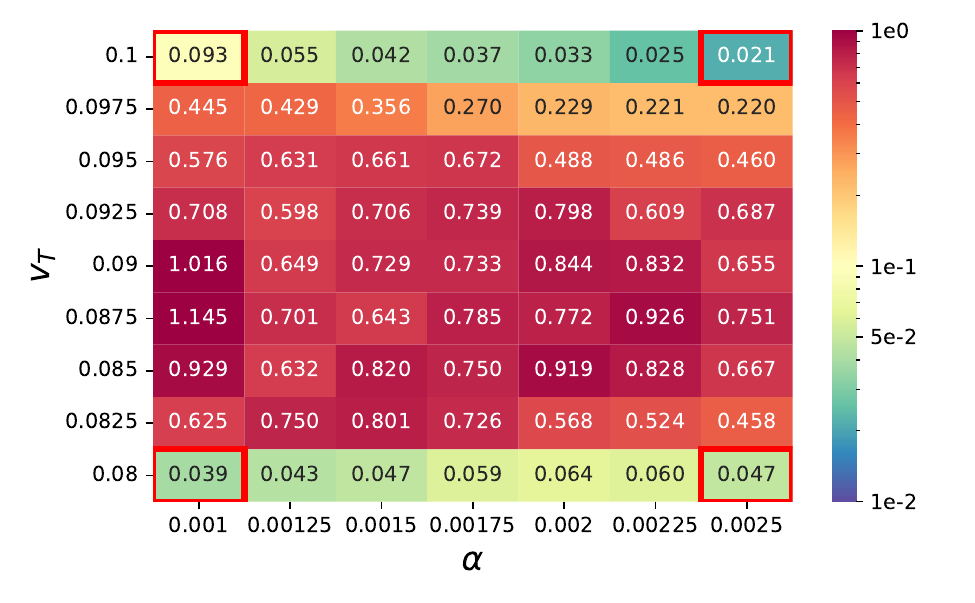}
    \end{subfigure}
    \begin{subfigure}{0.49\columnwidth}
    \caption{Single ROM, $\delta_\sigma = 5 \times \tento{-4}$}
    \includegraphics[width=1\columnwidth]{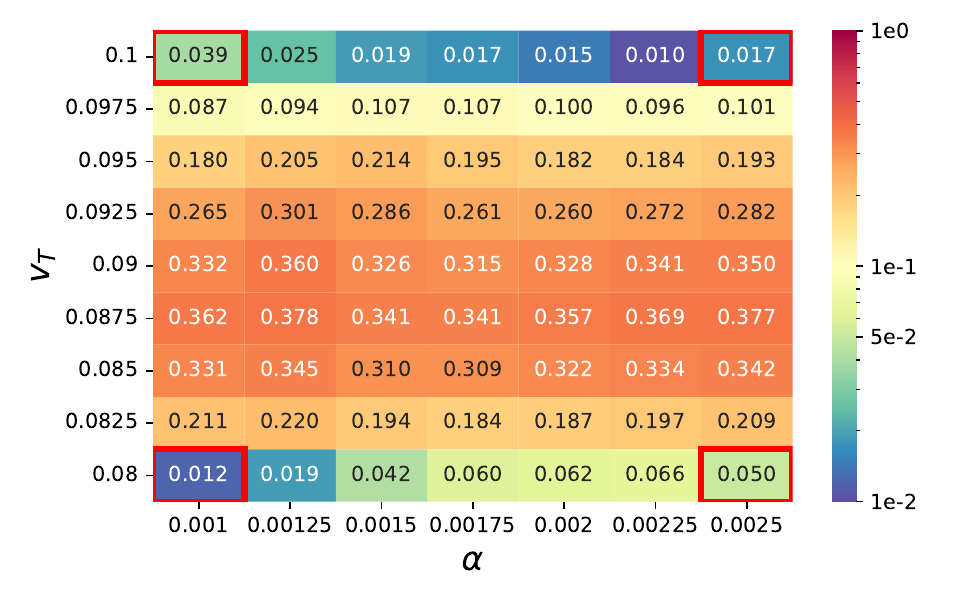}
    \end{subfigure}
    \begin{subfigure}{0.49\columnwidth}
    \caption{TW-ROM, $N_w=5$, $\delta_\sigma=\tento{-3}$}
    \includegraphics[width=1\columnwidth]{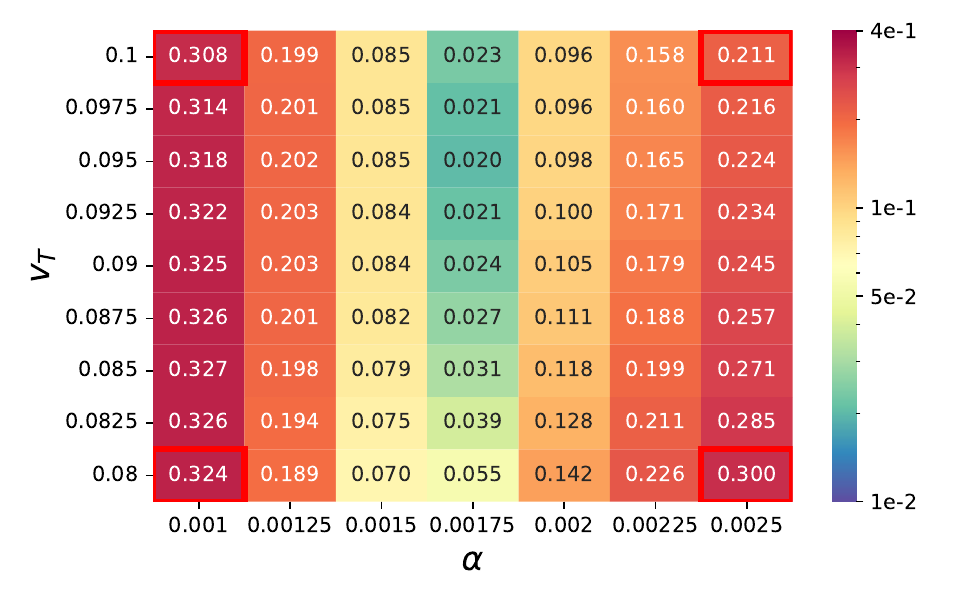}
    \end{subfigure}
    \begin{subfigure}{0.49\columnwidth}
    \caption{TW-ROM, $N_w=5$, $\delta_\sigma=5 \times \tento{-4}$}
    \includegraphics[width=1\columnwidth]{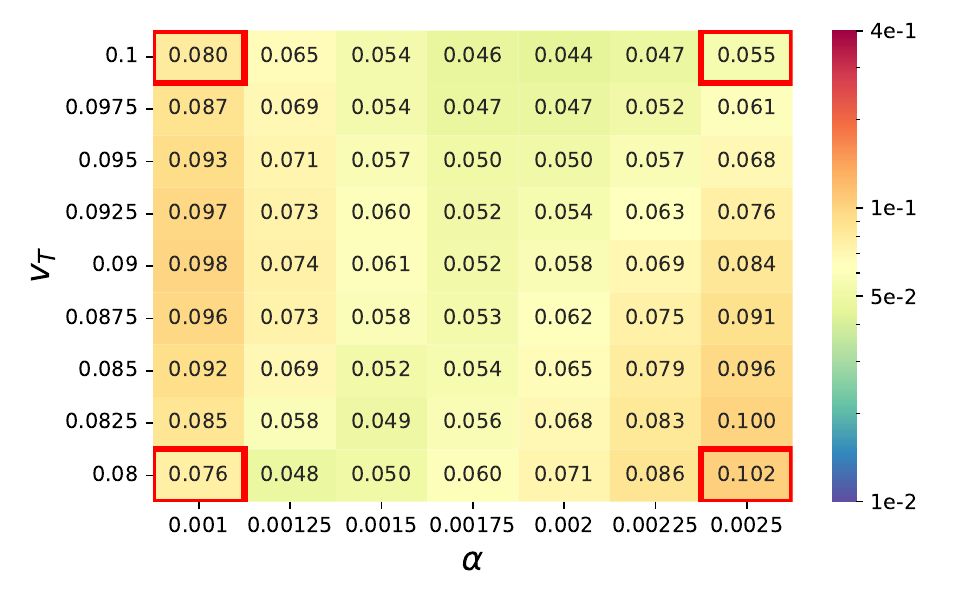}
    \end{subfigure}
    \begin{subfigure}{0.49\columnwidth}
    \caption{EW-ROM, $\delta_\sigma = \tento{-3}$}
    \includegraphics[width=1\columnwidth]{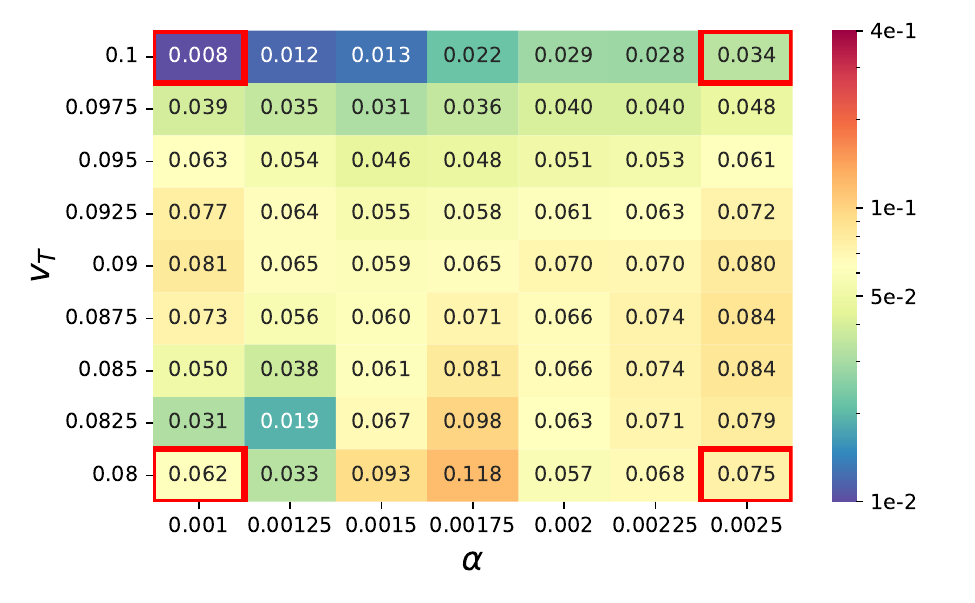}
    \end{subfigure}
    \caption{Parametric problem of the two-stream instability. The relative error in the distribution field at the final time $\epf$ of $\bp \in \mD^h$ for the single ROM, the TW-ROM with $N_w=5$ and the EW-ROM. 
    values of $\epf$ for the training parameters are enclosed with a red box.}
    \label{fig:pmor_4_compare_nw_4}
\end{figure}
In the online phase, the ROMs are deployed at the various parameters $\bp \in
\mD^h \subset \mD = [0.001,0.0025] \times [0.08,0.1]$ with final time $t_f = 10$,
where $\mD^h$ is a discrete parameter space $\mD^h \subset \mD$ %=
%[0.01,0.03]\times[0.8,1]$ 
with $9$ and $7$ evenly distributed discrete points in
the respective parameter range, resulting in a total of $63$ parameters. We note that the end-points $T_j(\bp)$ of the temporal subinterval $\mT_j(\bp)$ for both the TW-ROM and the EW-ROM are determined similarly as in the reproduction problem, which is described in the previous section.

Fig.~\ref{fig:pmor_4_compare_nw_4} displays the relative error in the distribution
field at the final time $\epf$ of $63$ parameters $\bp \in \mathcal{D}^h$ for
the single ROM, the TW-ROM, and the EW-ROM.
The results show that with $\delta_\sigma = \tento{-3}$, the EW-ROM achieves maximum errors of $7.5\%$ and $11\%$ in the training and testing parameters, respectively. On the other hand, although the single ROM is able to achieve a maximum error of $9\%$ in the training parameters, it performs terribly on almost all the testing parameters. The TW-ROM with $N_w=5$ has smaller errors in the testing parameters compared to the single ROM. However, surprisingly, it has the largest errors in the training parameters compared to the single ROM and the EW-ROM. We further decrease $\delta_\sigma$ from $\tento{-3}$ to $5 \times \tento{-4}$ for both the single ROM and the TW-ROM, and we find improvement in both ROMs. For the single ROM, the maximum error in the training parameters is $5\%$, and the maximum error in the testing parameters is reduced from $115\%$ to $38\%$. For the TW-ROM, the maximum error in the training parameters is $10\%$, and the maximum error in the testing parameters is reduced from $33\%$ to $10\%$.
We further calculate the average speed-up for the $63$ parameters in $\mathcal{D}^h$ for each ROM. We find that the single ROM with $\delta_\sigma= 5 \times \tento{-4}$ only achieved a speed-up of $16$. In contrast, we find that the TW-ROM with $\delta_\sigma = 5 \times \tento{-4}$ achieved a speed-up of $57$, while the EW-ROM achieved a speed-up of $89$.
Overall, considering both the accuracy and speed-up, the EW-ROM is a better model compared to the TW-ROM.

\begin{figure}[!ht]
    \centering
    \begin{subfigure}{0.49\columnwidth}
    \caption{$\bp = (0.001,0.1)$}
    \includegraphics[width=1\columnwidth]{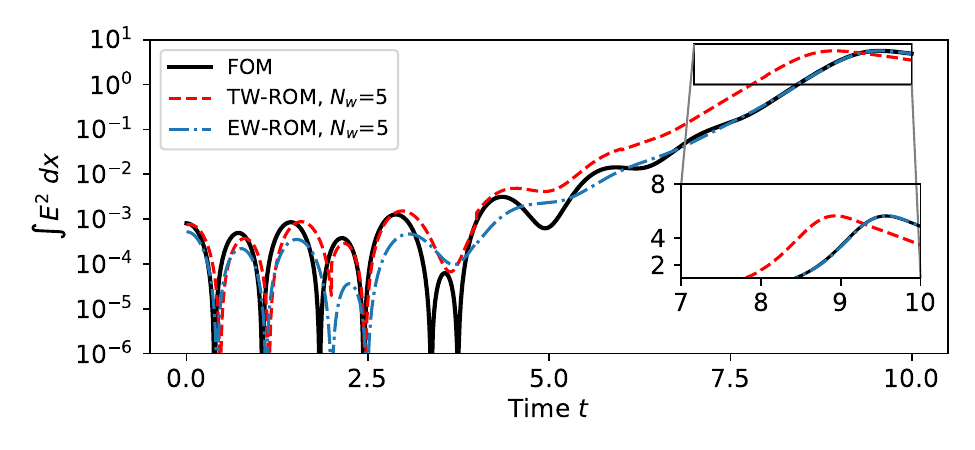}
    \end{subfigure}
    \begin{subfigure}{0.49\columnwidth}
    \caption{$\bp = (0.0025,0.08)$}
    \includegraphics[width=1\columnwidth]{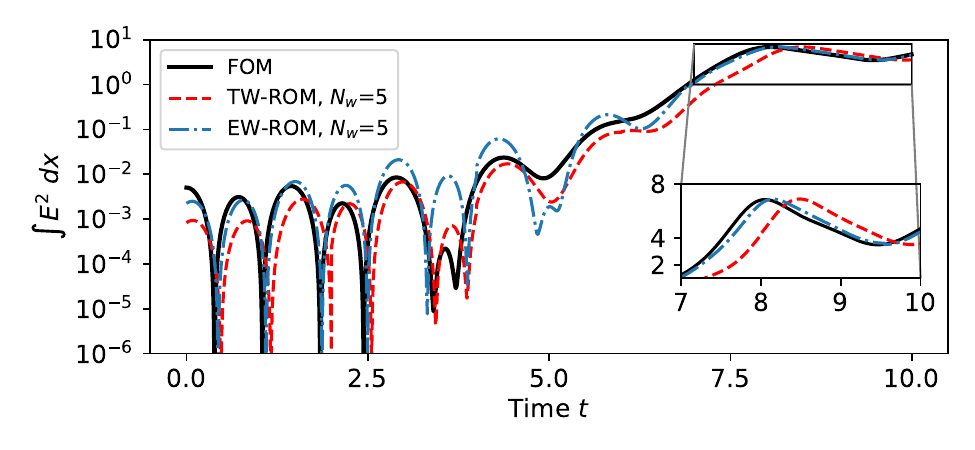}
    \end{subfigure}
    \begin{subfigure}{0.49\columnwidth}
    \caption{$\bp = (0.00175,0.08)$}
    \includegraphics[width=1\columnwidth]{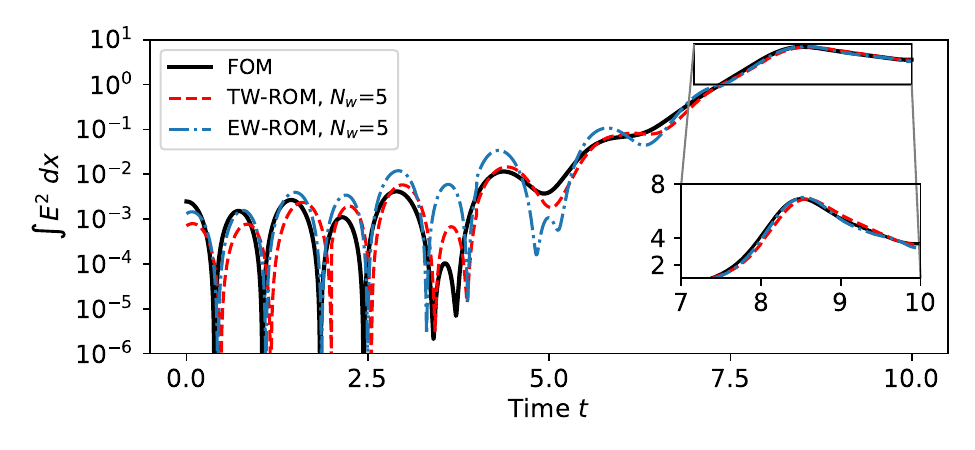}
    \end{subfigure}
    \begin{subfigure}{0.49\columnwidth}
    \caption{$\bp = (0.00175,0.1)$}
    \includegraphics[width=1\columnwidth]{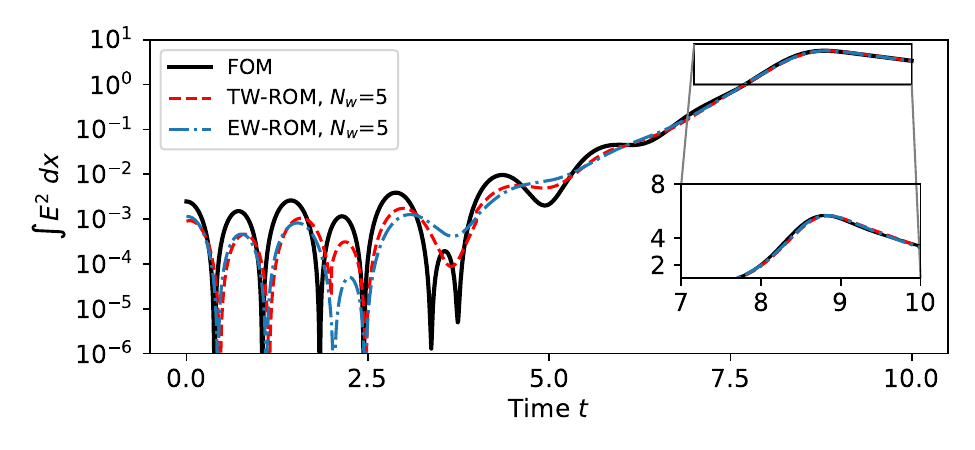}
    \end{subfigure}
    \begin{subfigure}{0.49\columnwidth}
    \caption{$\bp = (0.0025,0.075)$}
    \includegraphics[width=1\columnwidth]{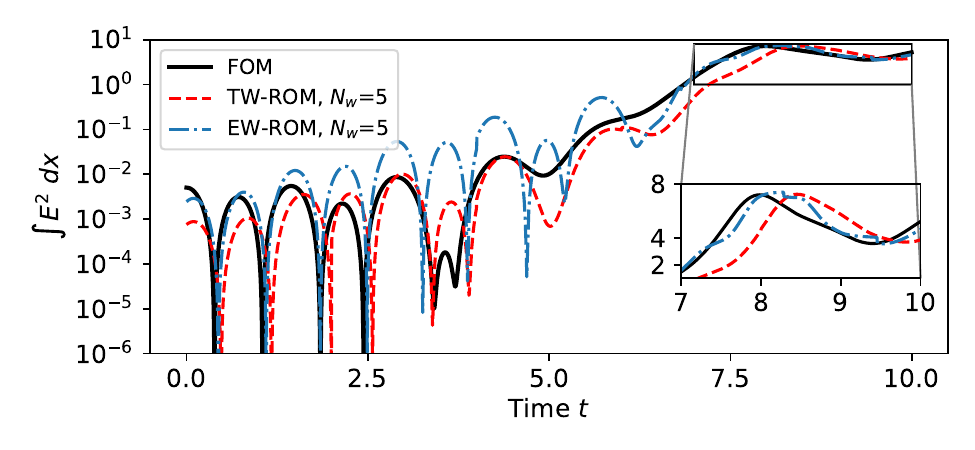}
    \end{subfigure}
    \begin{subfigure}{0.49\columnwidth}
    \caption{$\bp = (0.0025,0.0775)$}
    \includegraphics[width=1\columnwidth]{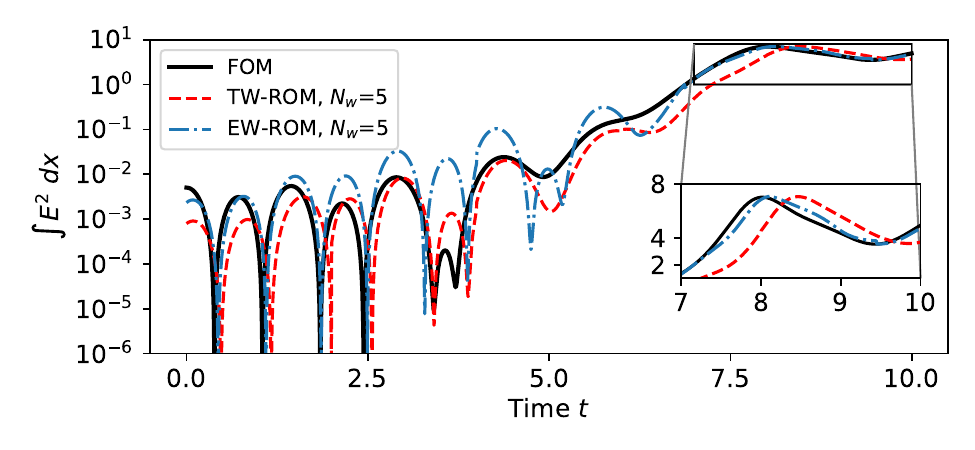}
    \end{subfigure}
    \caption{
    Comparison of the electric field energy $\int E^2~dx$ history between the EW-ROM and the TW-ROM with $N_w=5$ along with the FOM for the two training parameters $\bp = (0.001,0.1)$, and $(0.0025,0.08)$ and the four testing parameters $\bp=( 0.00175,0.08)$, $(0.00175,0.1)$, $(0.0025,0.075)$, and $(0.0025,0.0775)$. The energy missing ratio $\delta_\sigma$ is fixed to $\tento{-3}$.} 
    \label{fig:sc_pmor_e_compare}
\end{figure}

With $\delta_\sigma = \tento{-3}$, we further examine the behavior of the electric field energy $\int E^2~dx$
in the time interval $\Tinterval$ for the TW-ROM and the EW-ROM at the two training parameters  
$\bp=(0.001,0.1)$, $(0.0025,0.08)$ (training), and the two testing parameters $\bp = (0.00175,0.08)$, $(0.00175,0.1)$ in the interpolatory region and the two testing parameters $\bp = (0.0025,0.075)$ and $(0.0025,0.0775)$ in the extrapolatory region in Fig.~\ref{fig:sc_pmor_e_compare}.
For the training parameters, we find both the ROMs capture the growth rate of the instability, but the TW-ROM cannot capture the end-time of the growth stage. For the testing parameters in the interpolatory region, we find both ROMs capture the growth rate and end-time of the growth stage. For the testing parameters in the extrapolatory region, we find that TW-ROM is slightly better at capturing the growth rate than EW-ROM. However, we find that the TW-ROM struggles to predict the end-time of the growth stage, while the EW-ROM performs satisfactorily.

\begin{figure}[!ht]
    \centering
    \begin{subfigure}{0.49\columnwidth}
    \caption{$\bp = (0.001,0.08)$}
    \includegraphics[width=1\columnwidth]{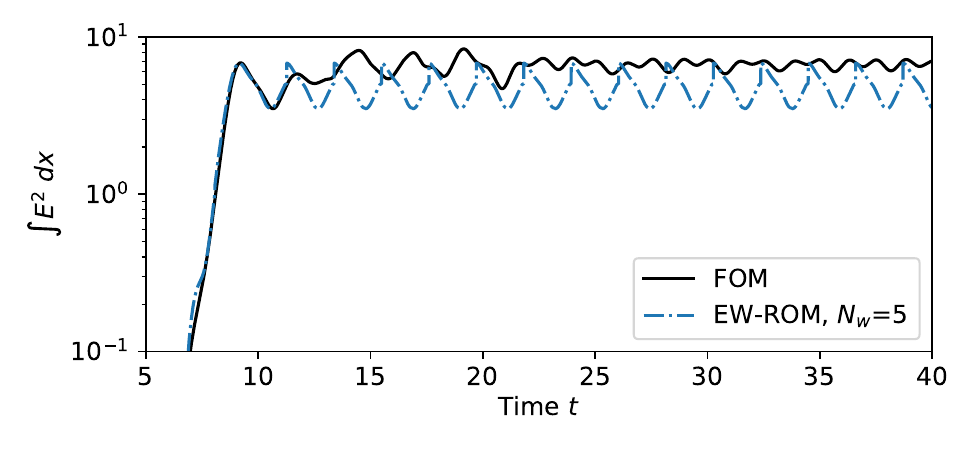}
    \end{subfigure}
    \begin{subfigure}{0.49\columnwidth}
    \caption{$\bp = (0.001,0.1)$}
    \includegraphics[width=1\columnwidth]{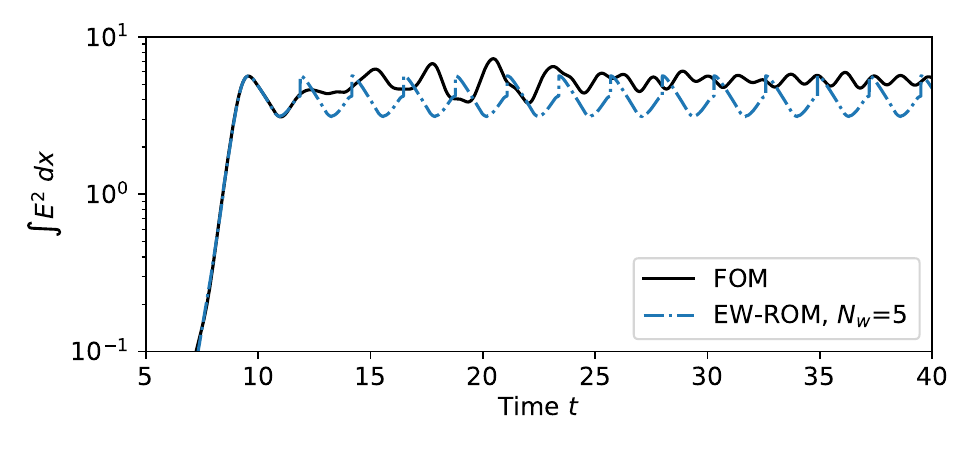}
    \end{subfigure}
    \begin{subfigure}{0.49\columnwidth}
    \caption{$\bp = (0.0025,0.08)$}
    \includegraphics[width=1\columnwidth]{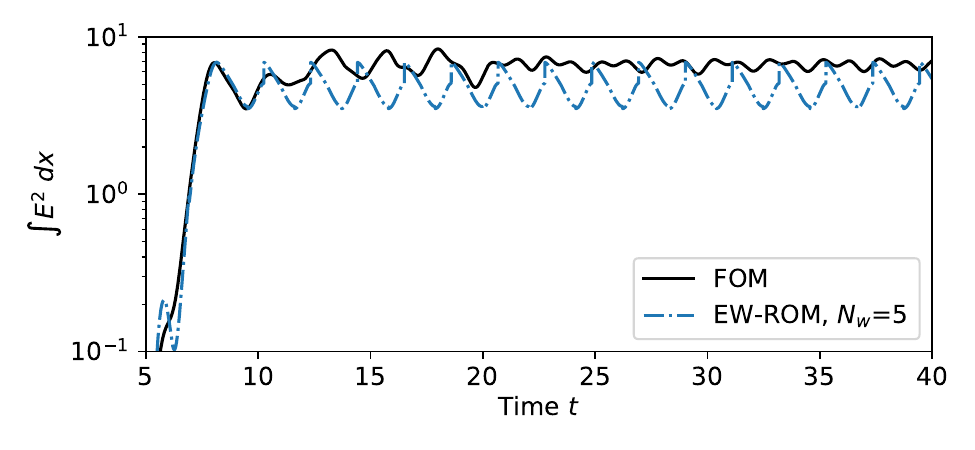}
    \end{subfigure}
    \begin{subfigure}{0.49\columnwidth}
    \caption{$\bp = (0.0025,0.1)$}
    \includegraphics[width=1\columnwidth]{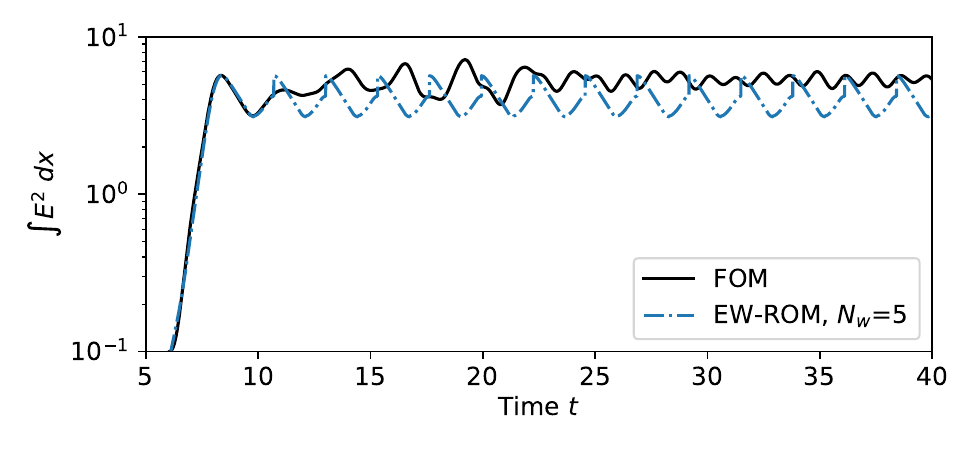}
    \end{subfigure}
    \caption{
    Extrapolation of the electric field energy $\int E^2~dx$ in time with the EW-ROM along with the FOM for the four training parameters. The energy missing ratio $\delta_\sigma$ is fixed to $\tento{-3}$.} 
    \label{fig:sc_pmor_e_ext}
\end{figure}

%In addition to using the ROM for predicting the electric field energy at the extrapolatory parameters $\bp$, 
We also investigate the performance of the ROM in extrapolation in time. We further simulate the EW-ROM for a time interval $[0, 40]$, which is four times larger than the training interval $\Tinterval=[0, 10]$. The electric field energy predictions made by the EW-ROM at four training parameters are shown in  Fig.~\ref{fig:sc_pmor_e_ext}. 
We find that, although EW-ROM does not accurately capture the energy, the predicted energy does a decent job of capturing the FOM energy accurately enough and does not blow up or decay. A more investigation would be required to further improve the performance, for example, by extending the training window.

\section{Conclusions}
\label{sec:conclusions}

In this paper, we develop reduced-order model techniques to accelerate the kinetic simulations of electrostatic plasmas. In particular, we study the effects of perturbation amplitude and thermal velocity as problem parameters on the physical quantities. 
We use POD to construct distribution and potential reduced basis functions.
Furthermore, we introduce a tensorial approach that  
utilizes the tensor-product representation in the grid-based simulation and reduced basis approximation of the electric field to efficiently update the nonlinear term. Moreover, we decompose the temporal domain into subintervals and construct a temporally local reduced-order model with varying parameters. We compare two decomposition techniques, namely, decomposing the solution manifold by the physical time and the electric field energy.
In Section~\ref{sec:numerical}, we demonstrate that the temporally-local ROMs are more accurate than the single ROM for the Landau-damping and the two-stream instability cases. In addition, in the two-stream instability problem, we show that the electric field energy indicator provides better results than the physical time indicator in terms of both speed-up and solution accuracy.

Possible directions of future investigation include applying the developed reduced-order model techniques to multiple spatial dimensions and extending them to the electromagnetic Vlasov\dash Maxwell equations.
We note that the solution dynamics in those contexts are expected to be more complex. This will require LS-ROMs to use a large number of modes to accurately represent the dynamic behavior.
Consequently, the development of a temporally-local ROM becomes essential. In addition, instead of constructing reduced spaces using POD, tensor decompositions could be beneficial in those high-dimensional systems \cite{hesthaven2023adaptive,guo2022low,hesthaven2022rank}. Furthermore, to extend the temporally local ROMs to a broader and higher-dimensional parameter space, the investigation of a more intelligent indicator is necessary.
The development of efficient error estimators to optimally select the training parameters to build ROMs is also important.

\section{CRediT Authorship Contribution Statement}

\textbf{Ping-Hsuan Tsai:}
Writing - original draft, Writing - review and editing,
Data curation, Investigation, Conceptualization, Validation, Visualization, Methodology, Software, Formal Analysis.
\textbf{Seung Whan Chung:}
Writing - review and editing, Conceptualization, Validation, Methodology, Supervision, Formal Analysis.
\textbf{Debojyoti Ghosh:} 
Writing - review and editing, Conceptualization, Methodology, Software, Supervision.
\textbf{John Loffeld:} 
Methodology, Software, Supervision.
\textbf{Youngsoo Choi:} 
Writing - review and editing, Conceptualization, Validation, Methodology, Supervision, Formal Analysis, Resources, Funding acquisition.
\textbf{Jonathan L. Belof:} 
Supervision, Project administration, Funding acquisition.

\section{Acknowledgment}
This work was performed under the auspices of the U.S. Department of Energy (DOE), by Lawrence Livermore National Laboratory (LLNL) under Contract No. DE-AC52–07NA27344 and was supported by Laboratory Directed Research and Development funding under projects 21-SI-006 and 22-SI-006.
Y.\ Choi was also supported for this work by the U.S. Department of Energy, Office of Science, Office of Advanced Scientific Computing Research, as part of the CHaRMNET Mathematical Multifaceted Integrated Capability Center (MMICC) program, under Award Number DE-SC0023164. IM release: LLNL-JRNL-871836.

\bibliographystyle{unsrt}
\bibliography{refs}

\begin{thebibliography}{10}

\bibitem{birdsall2018plasma}
Charles~K. Birdsall and A.~Bruce Langdon.
\newblock {\em Plasma physics via computer simulation}.
\newblock CRC press, 2018.

\bibitem{verboncoeur2005particle}
John~P. Verboncoeur.
\newblock Particle simulation of plasmas: review and advances.
\newblock {\em Plasma Physics and Controlled Fusion}, 47(5A):A231, 2005.

\bibitem{hypar}
Debojyoti Ghosh.
\newblock {HyPar} {R}epository, 2015.
\newblock {\tt https://github.com/debog/hypar}.

\bibitem{dorf2012progress}
MA~Dorf, RH~Cohen, JC~Compton, M~Dorr, TD~Rognlien, J~Angus, S~Krasheninnikov, P~Colella, D~Martin, and P~McCorquodale.
\newblock Progress with the {COGENT} edge kinetic code: collision operator options.
\newblock {\em Contributions to Plasma Physics}, 52(5-6):518--522, 2012.

\bibitem{lesser2022loki}
Steve Lesser, Alexey Stomakhin, Gilles Daviet, Joel Wretborn, John Edholm, Noh-Hoon Lee, Eston Schweickart, Xiao Zhai, Sean Flynn, and Andrew Moffat.
\newblock Loki: a unified multiphysics simulation framework for production.
\newblock {\em ACM Transactions on Graphics (TOG)}, 41(4):1--20, 2022.

\bibitem{vogman2014dory}
Genia~V. Vogman, Phillip Colella, and Uri Shumlak.
\newblock {D}ory--{G}uest--{H}arris instability as a benchmark for continuum kinetic {V}lasov--{P}oisson simulations of magnetized plasmas.
\newblock {\em Journal of Computational Physics}, 277:101--120, 2014.

\bibitem{hittinger2013block}
Jeffrey A.~F. Hittinger and Jeffrey~W. Banks.
\newblock Block-structured adaptive mesh refinement algorithms for {V}lasov simulation.
\newblock {\em Journal of Computational Physics}, 241:118--140, 2013.

\bibitem{sonnendrucker1999semi}
Eric Sonnendr{\"u}cker, Jean Roche, Pierre Bertrand, and Alain Ghizzo.
\newblock The semi-{L}agrangian method for the numerical resolution of the {V}lasov equation.
\newblock {\em Journal of Computational Physics}, 149(2):201--220, 1999.

\bibitem{einkemmer2016high}
Lukas Einkemmer.
\newblock High performance computing aspects of a dimension independent semi-{L}agrangian discontinuous {G}alerkin code.
\newblock {\em Computer Physics Communications}, 202:326--336, 2016.

\bibitem{rozza2022advanced}
Gianluigi Rozza, Giovanni Stabile, and Francesco Ballarin.
\newblock {\em Advanced reduced order methods and applications in computational fluid dynamics}.
\newblock SIAM, 2022.

\bibitem{fries2022lasdi}
William~D. Fries, Xiaolong He, and Youngsoo Choi.
\newblock {LaSDI}: Parametric latent space dynamics identification.
\newblock {\em Computer Methods in Applied Mechanics and Engineering}, 399:115436, 2022.

\bibitem{he2023glasdi}
Xiaolong He, Youngsoo Choi, William~D. Fries, Jonathan~L. Belof, and Jiun-Shyan Chen.
\newblock {gLaSDI}: Parametric physics-informed greedy latent space dynamics identification.
\newblock {\em Journal of Computational Physics}, 489:112267, 2023.

\bibitem{bonneville2024gplasdi}
Christophe Bonneville, Youngsoo Choi, Debojyoti Ghosh, and Jonathan~L. Belof.
\newblock {GPLaSDI}: Gaussian process-based interpretable latent space dynamics identification through deep autoencoder.
\newblock {\em Computer Methods in Applied Mechanics and Engineering}, 418:116535, 2024.

\bibitem{tran2024weak}
April Tran, Xiaolong He, Daniel~A Messenger, Youngsoo Choi, and David~M. Bortz.
\newblock Weak-form latent space dynamics identification.
\newblock {\em Computer Methods in Applied Mechanics and Engineering}, 427:116998, 2024.

\bibitem{bonneville2024comprehensive}
Christophe Bonneville, Xiaolong He, April Tran, Jun~Sur Park, William Fries, Daniel~A. Messenger, Siu~Wun Cheung, Yeonjong Shin, David~M. Bortz, Debojyoti Ghosh, et~al.
\newblock A comprehensive review of latent space dynamics identification algorithms for intrusive and non-intrusive reduced-order-modeling.
\newblock {\em arXiv preprint arXiv:2403.10748}, 2024.

\bibitem{kim2024gappy}
Youngkyu Kim, Youngsoo Choi, and Byounghyun Yoo.
\newblock {Gappy AE}: A nonlinear approach for gappy data reconstruction using auto-encoder.
\newblock {\em Computer Methods in Applied Mechanics and Engineering}, 426:116978, 2024.

\bibitem{choi2019space}
Youngsoo Choi and Kevin Carlberg.
\newblock Space--time least-squares {P}etrov--{G}alerkin projection for nonlinear model reduction.
\newblock {\em SIAM Journal on Scientific Computing}, 41(1):A26--A58, 2019.

\bibitem{carlberg2018conservative}
Kevin Carlberg, Youngsoo Choi, and Syuzanna Sargsyan.
\newblock Conservative model reduction for finite-volume models.
\newblock {\em Journal of Computational Physics}, 371:280--314, 2018.

\bibitem{kim2022fast}
Youngkyu Kim, Youngsoo Choi, David Widemann, and Tarek Zohdi.
\newblock A fast and accurate physics-informed neural network reduced order model with shallow masked autoencoder.
\newblock {\em Journal of Computational Physics}, 451:110841, 2022.

\bibitem{lauzon2024s}
Jessica~T. Lauzon, Siu~Wun Cheung, Yeonjong Shin, Youngsoo Choi, Dylan~M. Copeland, and Kevin Huynh.
\newblock S-{OPT}: A points selection algorithm for hyper-reduction in reduced order models.
\newblock {\em SIAM Journal on Scientific Computing}, 46(4):B474--B501, 2024.

\bibitem{diaz2024fast}
Alejandro~N. Diaz, Youngsoo Choi, and Matthias Heinkenschloss.
\newblock A fast and accurate domain decomposition nonlinear manifold reduced order model.
\newblock {\em Computer Methods in Applied Mechanics and Engineering}, 425:116943, 2024.

\bibitem{kaneko2020towards}
Kento Kaneko, Ping-Hsuan Tsai, and Paul Fischer.
\newblock Towards model order reduction for fluid-thermal analysis.
\newblock {\em Nuclear Engineering and Design}, 370:110866, 2020.

\bibitem{tsai2022parametric}
Ping-Hsuan Tsai and Paul Fischer.
\newblock Parametric model-order-reduction development for unsteady convection.
\newblock {\em Frontiers in Physics}, 10:903169, 2022.

\bibitem{hoang2021domain}
Chi Hoang, Youngsoo Choi, and Kevin Carlberg.
\newblock Domain-decomposition least-squares {P}etrov--{G}alerkin ({DD-LSPG}) nonlinear model reduction.
\newblock {\em Computer methods in applied mechanics and engineering}, 384:113997, 2021.

\bibitem{choi2021space}
Youngsoo Choi, Peter Brown, William Arrighi, Robert Anderson, and Kevin Huynh.
\newblock Space--time reduced order model for large-scale linear dynamical systems with application to {B}oltzmann transport problems.
\newblock {\em Journal of Computational Physics}, 424:109845, 2021.

\bibitem{choi2020gradient}
Youngsoo Choi, Gabriele Boncoraglio, Spenser Anderson, David Amsallem, and Charbel Farhat.
\newblock Gradient-based constrained optimization using a database of linear reduced-order models.
\newblock {\em Journal of Computational Physics}, 423:109787, 2020.

\bibitem{mcbane2021component}
Sean McBane and Youngsoo Choi.
\newblock Component-wise reduced order model lattice-type structure design.
\newblock {\em Computer methods in applied mechanics and engineering}, 381:113813, 2021.

\bibitem{mcbane2022stress}
Sean McBane, Youngsoo Choi, and Karen Willcox.
\newblock Stress-constrained topology optimization of lattice-like structures using component-wise reduced order models.
\newblock {\em Computer Methods in Applied Mechanics and Engineering}, 400:115525, 2022.

\bibitem{choi2019accelerating}
Youngsoo Choi, Geoffrey Oxberry, Daniel White, and Trenton Kirchdoerfer.
\newblock Accelerating design optimization using reduced order models.
\newblock {\em arXiv preprint arXiv:1909.11320}, 2019.

\bibitem{tsai2023accelerating}
Ping-Hsuan Tsai, Paul Fischer, and Edgar Solomonik.
\newblock Accelerating the {G}alerkin reduced-order model with the tensor decomposition for turbulent flows.
\newblock {\em arXiv preprint arXiv:2311.03694}, 2023.

\bibitem{tsai2024time}
Ping-Hsuan Tsai, Paul Fischer, and Traian Iliescu.
\newblock A time-relaxation reduced order model for the turbulent channel flow.
\newblock {\em Journal of Computational Physics}, page 113563, 2024.

\bibitem{zanardi2024scalable}
Ivan Zanardi, Alejandro~N. Diaz, Seung~Whan Chung, Marco Panesi, and Youngsoo Choi.
\newblock Scalable nonlinear manifold reduced order model for dynamical systems.
\newblock {\em arXiv preprint arXiv:2412.00507}, 2024.

\bibitem{cheung2023local}
Siu~Wun Cheung, Youngsoo Choi, Dylan~Matthew Copeland, and Kevin Huynh.
\newblock Local {L}agrangian reduced-order modeling for the {R}ayleigh-{T}aylor instability by solution manifold decomposition.
\newblock {\em Journal of Computational Physics}, 472:111655, 2023.

\bibitem{RAZAVI2025113576}
Alireza~H. Razavi and Masayuki Yano.
\newblock Registration-based nonlinear model reduction of parametrized aerodynamics problems with applications to transonic euler and rans flows.
\newblock {\em Journal of Computational Physics}, 521:113576, 2025.

\bibitem{washabaugh2012nonlinear}
Kyle Washabaugh, David Amsallem, Matthew Zahr, and Charbel Farhat.
\newblock Nonlinear model reduction for {CFD} problems using local reduced-order bases.
\newblock In {\em 42nd AIAA Fluid Dynamics Conference and Exhibit}, page 2686.

\bibitem{amsallem2012nonlinear}
David Amsallem, Matthew~J Zahr, and Charbel Farhat.
\newblock Nonlinear model order reduction based on local reduced-order bases.
\newblock {\em International Journal for Numerical Methods in Engineering}, 92(10):891--916, 2012.

\bibitem{copeland2022reduced}
Dylan~Matthew Copeland, Siu~Wun Cheung, Kevin Huynh, and Youngsoo Choi.
\newblock Reduced order models for {L}agrangian hydrodynamics.
\newblock {\em Computer Methods in Applied Mechanics and Engineering}, 388:114259, 2022.

\bibitem{cheung2024data}
Siu~Wun Cheung, Youngsoo Choi, H.~Keo Springer, and Teeratorn Kadeethum.
\newblock Data-scarce surrogate modeling of shock-induced pore collapse process.
\newblock {\em Shock Waves}, 34(3):237--256, 2024.

\bibitem{shimizu2021windowed}
Yukiko~S. Shimizu and Eric~J. Parish.
\newblock Windowed space--time least-squares petrov--galerkin model order reduction for nonlinear dynamical systems.
\newblock {\em Computer Methods in Applied Mechanics and Engineering}, 386:114050, 2021.

\bibitem{parish2021windowed}
Eric~J. Parish and Kevin~T. Carlberg.
\newblock Windowed least-squares model reduction for dynamical systems.
\newblock {\em Journal of Computational Physics}, 426:109939, 2021.

\bibitem{hesthaven2023adaptive}
Jan Hesthaven, Cecilia Pagliantini, and Nicol{\`o} Ripamonti.
\newblock Adaptive symplectic model order reduction of parametric particle-based {V}lasov--{P}oisson equation.
\newblock {\em Mathematics of Computation}, 2023.

\bibitem{jiangshu1996}
Guang-Shan Jiang and Chi-Wang Shu.
\newblock Efficient implementation of weighted {ENO} schemes.
\newblock {\em Journal of Computational Physics}, 126(1):202--228, 1996.

\bibitem{FFTW05}
Matteo Frigo and Steven~G. Johnson.
\newblock The design and implementation of {FFTW3}.
\newblock {\em Proceedings of the IEEE}, 93(2):216--231, 2005.
\newblock Special issue on ``Program Generation, Optimization, and Platform Adaptation''.

\bibitem{berkooz1993proper}
Gal Berkooz, Philip Holmes, and John~L. Lumley.
\newblock The proper orthogonal decomposition in the analysis of turbulent flows.
\newblock {\em Annual review of fluid mechanics}, 25(1):539--575, 1993.

\bibitem{fick2018stabilized}
Lambert Fick, Yvon Maday, Anthony~T. Patera, and Tommaso Taddei.
\newblock A stabilized pod model for turbulent flows over a range of reynolds numbers: Optimal parameter sampling and constrained projection.
\newblock {\em Journal of Computational Physics}, 371:214--243, 2018.

\bibitem{finn2023numerical}
Daniel~S. Finn, Matthew~G. Knepley, Joseph~V. Pusztay, and Mark~F. Adams.
\newblock A numerical study of {L}andau damping with {PETSc}-{PIC}.
\newblock {\em arXiv preprint arXiv:2303.12620}, 2023.

\bibitem{anderson2001tutorial}
D~Anderson, Renato Fedele, and M.~Lisak.
\newblock A tutorial presentation of the two stream instability and {L}andau damping.
\newblock {\em American Journal of Physics}, 69(12):1262--1266, 2001.

\bibitem{guo2022low}
Wei Guo and Jing-Mei Qiu.
\newblock A low rank tensor representation of linear transport and nonlinear {V}lasov solutions and their associated flow maps.
\newblock {\em Journal of Computational Physics}, 458:111089, 2022.

\bibitem{hesthaven2022rank}
Jan~S. Hesthaven, Cecilia Pagliantini, and Nicol{\`o} Ripamonti.
\newblock Rank-adaptive structure-preserving model order reduction of hamiltonian systems.
\newblock {\em ESAIM: Mathematical Modelling and Numerical Analysis}, 56(2):617--650, 2022.

\end{thebibliography}

\end{document}